\documentclass[a4paper,11pt]{amsart}
\usepackage{amsmath,amstext, amsthm, amssymb}
\usepackage[all]{xy}

\makeatletter
\numberwithin{equation}{section}

\usepackage[colorlinks=true, linkcolor=blue, 
hyperfootnotes=true,citecolor=blue,urlcolor=black]{hyperref}
 \usepackage{color}
 \usepackage{tikz,epic}
 \usepackage[all]{xy}
\usepackage{geometry,pgfplots}
\usepackage[english]{babel}
\reversemarginpar

\font \tenmat = msbm10 \font \sevenmat = msbm7 \font \fivemat = msbm5 \newfam \matfam \textfont \matfam = \tenmat \scriptfont
\matfam = \sevenmat \scriptscriptfont \matfam = \fivemat 

 \def\ZZ{{\mathbb Z}}  \def\RR{{\mathbb R}} \def\CC{{\mathbb C}} \def\HH{{\mathbb H}} \def\EE{{\mathbb E}}      
\def\DD{{\mathfrak D}}

\def\cB{{\mathcal B}} \def\cL{{\mathcal L}}
 
\def\cS{{\mathcal S}}

\def\sq{\sigma_{q}}

\newtheorem{thm}{Theorem}[section]
\newtheorem{lem}[thm]{Lemma}

\newtheorem{cor}[thm]{Corollary}
\newtheorem{prop}[thm]{Proposition}
\newtheorem{hypo}[thm]{Hypothesis}

\newtheorem{dfn}[thm]{Definition}
\newtheorem{rmk}[thm]{Remark}

\title{On the product of two generic $q$-Gevrey series}
\author{Thomas Dreyfus}
\address{Thomas Dreyfus, Université Bourgogne Europe, CNRS, IMB UMR 5584, F-21000 Dijon, France
}
\email{thomas.dreyfus@cnrs.fr}

\author{Changgui Zhang}
\address{Changgui Zhang, Laboratoire P. Painlev\'e CNRS UMR 8524, D\'epartement de Math\'ematiques, FST, Universit\'e de Lille, Cit\'e Scientifique, F-59655, Villeneuve d’Ascq cedex}
\email{changgui.zhang@univ-lille.fr}

\begin{document}
\date{\today}
\thanks{This project has received funding from the ANR de rerum natura ANR-19-CE40-0018 and also is supported by the Labex CEMPI (ANR-11-LABX-0007). The IMB receives support from the EIPHI Graduate School (contract ANR-17-EURE-0002).
}

\begin{abstract} 
In this paper, we consider a $q$-analog of the Borel-Laplace summation process introduced by Fabienne Marotte and the second author. We specifically examine two power series solutions of linear $q$-difference equations whose Newton polygon admits only positive slopes equal to $1$. These series, known as the generic $q$-Gevrey series, are shown to have the property that the product of two such series is $Gq$-summable at double level $(1,2)$. Furthermore, we prove that the $Gq$-sum of this product equals the product of the $Gq$-sums of the original two series.
\end{abstract}

\maketitle

\tableofcontents

\section{Introduction}\label{sec:primo}
This paper deals with the algebraic properties of a $q$-analog of the Borel-Laplace summation process defined in \cite{marotte2000multisommabilite}. Before considering the $q$-world, let us take a quick look at the theory in the setting  of linear differential equations. We refer to \cite{balser2006divergent, balser2008formal} for a detailed explanation of this theory.\par 

Consider a meromorphic linear differential equation. We have the coexistence of divergent formal power series and integral solutions. For example, the Euler equation $x^{2}\partial_x y+y=1$  admits the Euler series $ \hat E (x) :=\displaystyle \sum_{n=0}^{\infty} (-1)^{n}n! x^{n}$ as formal power series solution. On the other side, there are integral solutions like $$\mathcal{S}^{d}( \hat E )(x):=\displaystyle\int_{0}^{\infty e^{\mathbf{i}d}}\frac{e^{-\xi/x}}{1+\xi}d\xi .$$

The integration path is an half-line in $\CC$ with complex numbers having an argument $d\in \RR$. This integral is well-behaved when the integration path avoids the pole $\xi=-1$, meaning $d\not \equiv \pi [2\pi]$. It can be proven that the function is analytic in the sector $\arg(x)\in (d-\pi/2,d+\pi/2)$ and is asymptotic to $\hat{E}$ in the Gevrey sense.

More generally, for any formal power series solution of a linear differential equation whose coefficients are germs of meromorphic functions at $0$, we can, for a convenient $d\in \RR$, construct an integral solution using Borel and Laplace transformations in the direction of $d$. The map $f\mapsto \mathcal{S}^{d}(f)$, which transforms a formal power series into the integral solution, acts as a morphism of fields. Moreover, it preserves the germs of meromorphic functions at $0$ and commutes  with derivation, meaning $\mathcal{S}^{d}(\partial_x f)=\partial_x \mathcal{S}^{d}(f)$. \\ \par 

Now, let us shift our focus to $q$-difference equations. Through the whole paper, $q$ denotes a real number in the open interval $(1,+\infty)$.  Before going further, it is worthwhile to recollect some facts about one of the simplest generic $q$-Gevrey series, that will be seen as a $q$-analog of the well-known Euler series:
\begin{equation}\label{equation:qEuler}
 \hat E_q(x)=\sum_{n\ge 0}(-1)^nq^{n(n-1)/2}\,x^n\,.
\end{equation} 
Defining the $q$-difference operator as $\sq y(x):= y(qx)$, it is solution of the $q$-Euler equation $x\sq y+y=1$. A convergent  solution takes the form 
$$\mathcal{S}_{q}^{d}(\hat{E}_q):= \frac{1}{\sqrt{2\pi\log (q)}}\,\displaystyle\int_0^{\infty e^{\mathbf{i}d}}\frac{e^{-(\log(\frac{x}{\sqrt{q}\xi}))^2/2\log(q)}}{1+\xi}\,\,\frac{d\xi}{\xi}, \quad d\not \equiv \pi [2\pi].$$
It is worth mentioning that this solution is meromorphic on the Riemann surface of the logarithm and defines a multivalued complex function. \par 
 Besides, it is easy to see that $\hat E_q(x)^2$ satisfies a non-homogeneous second order $q$-difference equation; for instance, see \cite[\S 2.2.2]{marotte2000multisommabilite} or \cite[(1.3)]{zhang2023positive}. Thus, by invoking \cite[Theorem 3.3.5]{marotte2000multisommabilite}, it can be found that the square $\hat E_q(x)^2$ is $Gq$-multisummable,  that is, we can apply a $q$-Borel-Laplace summation process to it to obtain a meromorphic solution having the series as an asymptotic expansion in the $q$-Gevrey sense. Specifically $\hat E_q(x)^2$ is $Gq$-summable of double level $(1, 2)$; while $\hat E_q(x)$ is $Gq$-multisummable of level~$1$. Further details on this definition will be provided later. \par 

With this established, a pertinent inquiry arises: does the $Gq$-multisum of the square of the $q$-Euler series equate to the square of the $Gq$-sum of the initial series?
This paper will address this question in a broader context as detailed later.\par 

Unfortunately, unlike the differential case, there is not a  unique $q$-summation process due to the nonuniqueness of the $q$-analog of the exponential. Multiple $q$-sums exist, each with its advantages in specific contexts. We can find more details in the works of the two authors, Ramis, Sauloy, and a more exhaustively \cite{jackson1910xxii,jackson1917q,abdi1960q,abdi1964certain,
ramis1992growth,zhang1999developpements,
marotte2000multisommabilite,zhang2000transformations,
zhang2001fonction,ramis2002developpement,
zhang2002sommation,
zhang2003fonctions,zhang2004solutions,di2009q,
ramis2013local,dreyfus2015building,
dreyfus2016q}. 
In this paper, we use a $q$-summation process introduced in \cite{marotte2000multisommabilite}. The main reason for this choice is that we are interested in the algebraic properties of the $Gq$-sum. Note that the question was already considered in \cite{jackson1910xxii,jackson1917q}, but the $q$-summation process considered in this paper seems to be the best candidates as shows the following results.
It has been shown that given a formal power series solution $f$ of a linear meromorphic $q$-difference equation, we can associate a multivalued meromorphic solution $\mathcal{S}_{q}^{d}(f)$ using $q$-analogues of the Borel-Laplace summation in the direction $d\in \RR$, where $d$ is chosen conveniently.

The authors established that the $q$-summation process $f\mapsto \mathcal{S}_{q}^{d}(f)$ satisfies certain algebraic properties, as outlined in \cite{marotte2000multisommabilite}.
\begin{itemize}
\item  If $f_1,f_2$ are $Gq$-summable in the direction $d\in \RR$, then   $f_1+f_2$ $Gq$-summable in the direction $d$ and $\mathcal{S}_{q}^{d}(f_{1}+ f_{2})=\mathcal{S}_{q}^{d}(f_{1})+ \mathcal{S}_{q}^{d}(f_{2})$;
\item  If $f$ is $Gq$-summable in the direction $d$,  then $\sigma_{q}(f)$ is  $Gq$-summable in the direction $d$ and 
$\mathcal{S}_{q}^{d}(\sigma_{q}(f))=\sigma_{q}\left(\mathcal{S}_{q}^{d}(f)\right)$;
\item For all convergent series $f$,   if $g$ is $Gq$-summable in the direction $d$,  then $fg$ is $Gq$-summable in the direction $d$ and   $\mathcal{S}_{q}^{d}(fg)=f\mathcal{S}_{q}^{d}(g)$.
\end{itemize}
 Let $\CC[[x]]$ be the $\CC$-vector space of power series of the variable $x$, and let $\CC\{x\}$ be the $\CC$-vector space of convergent power series, which can be identified with the set of germs of analytic functions at $x=0$ in the complex plane.  As in \cite[\S 4]{di2009q}, we call {\it a generic $q$-Gevrey series} a power series $\hat f\in\CC[[x]]$ that satisfies a linear $q$-difference equation $\Delta \hat f=g\in\CC\{x\}$ whose Newton polygon admits only one finite slope equal to one. More explicitly, one can and one will assume that $\Delta$ represents a certain $n$-th order $q$-difference operator of the following form : 
\begin{equation}\label{equation:Delta}
 \Delta=\sum_{k=0}^na_{n-k}(x)(x\,\sigma_q)^k\in\CC\{x\}[x\sigma_q],\quad a_0(0)\,a_n(0)\not=0\,,
\end{equation}
Recall that in virtue of \cite{marotte2000multisommabilite} every generic $q$-Gevrey series are $Gq$-summable.
In this paper, we establish that, under specific assumptions, the $q$-summation process and the product operation commute.  Further details on the meaning of $\mathcal{S}_{q;1}^{d},\mathcal{S}_{q;(1,2)}^{d}$ will be provided later but roughly speaking it means that $\mathcal{S}_{q;1}^{d}$ corresponds to the $Gq$-sum of a series that is $Gq$-summable of level $1$, and $ \mathcal{S}_{q;(1,2)}^{d}$ corresponds to the $Gq$-sum of a series that is $Gq$-summable of level $(1,2)$.

\begin{thm}\label{theorem:intro} 
Let $(\hat f_1,\hat f_{2})$ be a pair of generic $q$-Gevrey series, $Gq$-summable in the direction $d\in \RR$. Then the product  $\hat f_1\,\hat f_2$ is $Gq$-summable in the direction $d$. Furthermore, the following holds:
\begin{equation}
 \label{equation:f1f2}
 \mathcal{S}_{q;(1,2)}^{d}(\hat f_{1}\, \hat f_{2})=\mathcal{S}_{q;1}^{d}\,(\hat f_{1})\,\mathcal{S}_{q;1}^{d}(\hat f_{2}).
\end{equation}
\end{thm}
In this work, we have chosen to focus on the case where $q$ is real. Although an analogue of our results is expected to hold for non-real $q$, with $|q|>1$ would require introducing a generalized $q$-summation. In that setting, the integration contour would no longer be a half-line but rather a spiral of the form $q^{\mathbb{R}}e^{\mathbf{i}d}$, which would significantly complicate the convergence domains throughout the paper. This practical difficulty motivates our restriction to real $q$.
We  conjecture that $f\mapsto \mathcal{S}_{q}^{d}(f)$ defines a ring morphism, and this result is a step in that direction. Proving this conjecture could help us to define the $q$-analog of the Stokes operators in a manner analogous to the differential case. It may also lead to proving a $q$-analog of the Ramis density theorem.
It is worth mentioning that  the latter question has been considered from a different point of view in \cite{ramis2013local}, and a comparison between the two approaches would be particularly interesting.  Finally, we believe the present work may contribute to a better understanding of the analytical structures underlying certain non-linear problems, including $q$-Painlevé equations and conformal field theory. For more details, see \cite{jimbo2017cft}.\par 

The remainder of the paper is structured as follows:

{\bf Section \ref{section:revisit}:} We will review the concepts of $Gq$-summation and its variant of double-level extension as introduced in \cite{zhang1999developpements} and \cite{marotte2000multisommabilite}, respectively. Additionally, we will define the so-called generalized $q$-Euler series by differentiating  a variant of the $q$-Euler series with respect to a parameter. This allows us to state a new characterization for any generic $q$-Gevrey series; see Theorem \ref{theorem:generic}. 

{\bf Section \ref{section:qGevrey}:} We will first present several preliminary lemmas concerning the generalized $q$-Euler series. By revisiting a result from \cite{di2009q}, we will give a proof of Theorem~\ref{theorem:generic}. 

{\bf Section \ref{section:qEuler}:} This section is dedicated to establishing Theorem \ref{theorem:intro} for products of generalized $q$-Euler series; see Theorem~\ref{theorem:Eab}. For  achieve this, particular attention will be given for the study of such products in two successive $q$-Borel planes. 

{\bf Section \ref{section:proof}:} We will focus on proving the main result of the paper, namely Theorem~\ref{theorem:intro}. This will be done in a somewhat technical manner, specifically combining Theorems~\ref{theorem:generic} and \ref{theorem:Eab}.

\section{Revisiting $Gq$-summation and generic $q$-Gevrey series}\label{section:revisit}

 Throughout this section, we will adhere to the notation for functional transformations and spaces introduced in \cite{zhang1999developpements, marotte2000multisommabilite}, applying it wherever relevant. As usual, $\tilde\CC^*$ denotes the Riemann surface associated with the logarithm function, while $\log$ signifies the principal branch of this function. Specifically, for any real number \( a > 0 \), the notation $\log a$ corresponds to the natural logarithm $\ln a$.
 
\subsection{Space of $Gq$-summable series $\CC\{x\}_{q;1}$}\label{subsection:Gq1}  
Let $\CC[[x]]_{q;1}$ be a space of $q$-Gevrey series of the first order, {\it i.e.}:
\begin{equation}\label{equation:qGevry}
\sum_{n\ge 0}c_nx^n\in\CC[[x]]_{q;1}\quad\Longleftrightarrow\quad \exists\, C,A>0, \forall\, n\ge 0,\, |c_n|\le C\,A^nq^{n(n-1)/2}.
\end{equation}
 Thus, it is clear that the inclusion relation $\CC\{x\}\subset\CC[[x]]_{q;1}\subset\CC[[x]]$ holds.
 To simplify, we define
\begin{equation}\label{equation:e_q}
 e_q(x)=\sqrt{2\pi\log (q)}\,e^{\log(\sqrt{q}\,x)^2/(2\log(q))}\quad (x\in\tilde\CC^*)\,,
\end{equation}
and
\begin{equation}\label{equation:DV}
\tilde D_{R}=\{x\in\tilde\CC^*: |x|<R\}\,,\quad 
V_{\varepsilon}^d=\left\{\xi\in\tilde\CC^*: |\arg (\xi)-d|<\varepsilon\right\}\,,
\end{equation}
where $R>0$, $d\in\RR$ and $\varepsilon>0$. Here, $d$ represents the direction of the argument, and we will refer to $d$ simply as the direction.  In addition, we will denote by $\mathcal{C}^d$ the analytic continuation along the direction $d$ in some sector $V_\epsilon^d$. 
  
\begin{dfn}\label{definition:Gq1}
 A power series $\hat f= \sum_{n\ge 0}a_nx^n$ is said to be {\it $Gq$-summable in a given direction $d$} and one writes $\hat f\in\CC\{x\}_{q;1}^d$, if one has $\hat \cB_{q;1}\hat f\in\HH_{q;1}^d$, where $\hat \cB_{q;1}\hat f$ and  $\HH_{q;1}^d$ are defined as follows:
 \begin{itemize}
  \item 
 $\hat \cB_{q;1}\hat f(\xi)=\sum_{n\ge 0}a_n\,q^{-n(n-1)/2}\,\xi^n$, and this is what we call {\it the (first order) formal $q$-Borel transform of $\hat f$};
 \item $\HH_{q;1}^d$ denotes the $\CC$-vector space of germs of analytic functions at $\xi=0\in\CC$ that can be continued into an analytic function $\phi$ possessing a $q$-exponential growth of order one at the infinity in some sector $V_\epsilon^d$ in the following sense: 
 \begin{equation}
  \label{equation:qexp1}\exists A>0,\quad \phi(\xi)=O\left(e_q(A\xi)\right)\quad \hbox{as}\quad |\xi|\to\infty\quad \hbox{with } \quad \xi\in V_{\epsilon}^{d}\,.  
 \end{equation}
 \end{itemize}
 In this case, the $Gq$-sum of $\hat f$ related to the direction $d$ is, by definition, the $q$-Laplace transform of $\phi$ in the following manner: $\phi=\mathcal{C}^d\circ\hat\cB_{q;1}^d\hat f$, $\mathcal{S}_{q;1}^d\hat f=\mathcal{L}_{q;1}^d\phi$, and
 \begin{equation}
  \label{equation:qLaplace1}
  \mathcal{L}_{q;1}^d\phi(x)=\int_0^{\infty e^{id}}\frac{\phi(\xi)}{e_q({\xi}/{x})}\,\frac{d\xi}{\xi}\,,
 \end{equation}
where $x$ belongs to some domain $\tilde D_R$ ($R>0$) of the Riemann surface $\tilde\CC^*$.

Furthermore, if there exists a finite subset $S\subset [0,2\pi]$ such that $\hat f\in\CC\{x\}_{q;1}^d$ for any $d\in\left([0,2\pi]\setminus S\right)$, we will say that $\hat f$ is {\it $Gq$-summable} and simply write $\hat f\in\CC\{x\}_{q;1}$. 
\end{dfn}

In view of \cite[Proposition 4.3.1 \& Théorème 5.3]{zhang1999developpements}, one has the following result.
\begin{prop}\label{prop:Gqsummable} 
Let $d\in \mathbb{R}$. 
\begin{enumerate}
\item Given $\hat f\in\CC\{x\}_{q;1}^d$ and $h\in\CC\{x\}$, one has $h\,\hat f\in\CC\{x\}_{q;1}^d$ and $\cS_{q;1}^d(h\,\hat f)=h\,\cS_{q;1}^d\hat f$.
\item 
Any generic $q$-Gevrey power series is  $Gq$-summable.
\end{enumerate}
\end{prop}

\subsection{Space of $Gq$-summable series $\CC\{x\}_{q;(1,2)}$}\label{subsection:Gq2} 
In order to simplify the presentation, we will only recall the defnition of the double-level $Gq$-summable series ; for the more general definition, see \cite{marotte2000multisommabilite}.

\begin{dfn}\label{definition:Gq2}
 A power series $\hat f= \sum_{n\ge 0}a_nx^n$ is said to be {\it $Gq$-summable of double-level $(1,2)$ in a given direction $d$} and one writes $\hat f\in\CC\{x\}_{q;(1,2)}^d$, if   $\hat \cB_{q;1}\hat f\in\HH_{q;2}^d$ and $\cL_{q;2}^d\circ\mathcal{C}^d\circ\hat \cB_{q;1}\hat f\in\tilde\HH_{q;2}^d$, where  $ \HH_{q;2}^d$, $\cL_{q;2}^d$ and  $\tilde\HH_{q;2}^d$ are defined as follows:
 \begin{itemize}
  \item 
$\HH_{q;2}^d=\HH_{q^{1/2};1}^d$ and  $\cL_{q;2}^d=\cL_{q^{1/2};1}^d$, that are obtained by replacing $q$ with $q^{1/2}(=\sqrt q)$ in \eqref{equation:qexp1} and \eqref{equation:qLaplace1}, respectively;
 \item $\tilde\HH_{q;2}^d$ denotes the $\CC$-vector space of analytic functions $\phi$ in some sector $V_\epsilon^d$ ($\epsilon>0$ being arbitrary) such that $\phi(\xi)=O(1)$  as $\xi\to 0$ in $V_\epsilon^d$ and that the condition in \eqref{equation:qexp1} is fulfilled for $q^{1/2}$ instead of $q$. 
 \end{itemize}
 If this is the case, the corresponding  $Gq$-sum of $\hat f$, $\mathcal{S}_{q;(1,2)}^d\hat f$, is defined as follows: $\phi_0=\mathcal{C}^d\circ\hat\cB_{q;1}^d\hat f$, $\phi_1^d=\mathcal{C}^d\circ\cL_{q;2}^d\phi_0$, and $\mathcal{S}_{q;(1,2)}^d\hat f=\mathcal{L}_{q;2}^d\phi_1^d$.
 
Furthermore, if  $\hat f\in\CC\{x\}_{q;(1,2)}^d$ for any value of $d$ in $[0,2\pi]$ but eventually excepted one finite number, we will say that $\hat f$ is {\it $Gq$-summable of level $(1,2)$} and simply write $\hat f\in\CC\{x\}_{q;(1,2)}$. 
\end{dfn}

By considering \cite[Lemme 2.4.1 \& Théorème 2.4.3]{marotte2000multisommabilite}, one gets the following statement.

\begin{prop}\label{prop:2Gqsummable} Given any $d\in\RR$, the set
$\CC\{x\}_{q;(1,2)}^d$ constitutes a $\CC\{x\}$-module containing $\CC\{x\}_{q;1}^d$ as a sub-module.

\end{prop}

\subsection{Generic $q$-Gevrey series in terms of $q$-Euler series}
It is easy to see that every convergent power series is a generic $q$-Gevrey series. 
The $q$-Euler series  $\hat E_q$ given in \eqref{equation:qEuler} can be seen as one of the simplest generic $q$-Gevrey series whose radius of convergence is zero. Indeed, a directe calculation shows that 
\begin{equation}\label{equation:qEulerEquation}
 (x\,\sigma_q+1)\hat E_q(x)=1\,.
\end{equation}
To any $a\in\CC^*$ will be associated the family $(\hat E_q^{(a,\alpha)})_{\alpha\in\ZZ_{\ge 0}}$ of the so-called {\it generalized $q$-Euler series} as follows: $\displaystyle\hat E_q^{(a,0)}(x)=-\frac{1}{a}\,\hat E_q(-\frac{x}{a})$, and
\begin{equation}\label{equation:qEulera}
 \hat E_q^{(a,\alpha)}(x)=\frac{1}{\alpha !}\,{\partial _a^\alpha}\hat E_q^{(a,0)}(x)
\end{equation}
for $\alpha\ge 1$.
It follows that $\hat E_q=\hat E_q^{(-1,0)}$. Furthermore, for any $\lambda\in\CC^*$, one gets from considering \eqref{equation:qEulera} that
\begin{equation}\label{equation:qEuleraa'}
 \hat E_q^{(a/\lambda,\alpha)}(x)=\lambda^{\alpha+1}\,\hat E_q^{(a,\alpha)}(\lambda\,x)\,.
\end{equation}

\begin{thm}\label{theorem:generic} 
A power series $\hat f\in\CC[[x]]$ is a generic $q$-Gevrey series if, and only if, there exist $m\in\ZZ_{\ge0}$,  $\{f_0, f_1,\dots,f_m\}\subset\CC\{x\}$, $\{a_1, \dots,a_m\}\subset \CC^*$ and $\{\alpha_1,\dots,\alpha_m\}\subset\ZZ_{\ge0}$ such that the following equality holds in $\CC[[x]]$:
  \begin{equation}\label{equation:decomposition}
   \hat f=f_0+f_1\,\hat E_q^{(a_1,\alpha_1)}+\dots+f_m\,\hat E_q^{(a_m,\alpha_m)}\,.
  \end{equation}
\end{thm}

The proof of Theorem \ref{theorem:generic}, which is closely linked to \cite[Theorem 4.20]{di2009q}, will be given in \S \ref{subsection:thmGeneric}. Let's emphasize that this result plays a central role in the proof of Theorem \ref{theorem:intro}. Indeed, thanks to this result, it will suffice to verify the theorem for a pair of generalized $q$-Euler series.

\section{Generic $q$-Gevrey series in terms of the generalized $q$-Euler series}\label{section:qGevrey}

This section aims to prove Theorem \ref{theorem:generic} by first establishing properties for every family of generalized $q$-Euler series $\{\hat E_q^{(a/q^k,\alpha)}\}_{k,\alpha\in\ZZ_{\ge 0}}$. In \S \ref{subsection:qEuler1}, we will begin by examining the functional equation satisfied by each generalized $q$-Euler series \(\hat{E}_q^{(a, \alpha)}\), noting that it is $Gq$-summable in any direction that does contain the point with affix \(a\). Furthermore, it will be shown that every $Gq$-sum is bounded in any sector of finite openness; see Remark~3.2. 

In \S \ref{subsection:qEuler2}, we will consider the linear relationships that exist among any finite family \(\{ \hat{E}_q^{(a, \alpha)} \}\) where \(a \in a_0 q ^\mathbb{Z}\) and \(\alpha \in (\alpha_0 + \mathbb{Z})\). After that, in \S \ref{subsection:GenericOperators}, we will explain how to construct a $q$-difference operator associated with a given linear combination of generalized $q$-Euler series; see Theorem \ref{theorem:Eaqk}. In this way, we will provide a new characterization for a power series to be generic $q$-Gevrey, thereby proving Theorem~\ref{theorem:generic}.

\subsection{Functional equations and $Gq$-sums of the generalized $q$-Euler series}\label{subsection:qEuler1}

One recalls that the $q$-Euler series $\hat E_q$, which is precisely $\hat E_q^{(-1,0)}$, satisfies the $q$-difference equation given in \eqref{equation:qEulerEquation}. This result can be extended to any generalized $q$-Euler series $\hat E_q^{(a,\alpha)}$.

\begin{lem}\label{lemma:Ea}
For any $(a,\alpha)\in\CC^*\times\ZZ_{\ge 0}$, the generalized $q$-Euler series $\hat E_q^{(a,\alpha)}$ defined in \eqref{equation:qEulera} satisfies the $q$-difference equation
 \begin{equation}
 \label{equation:EaEq}
 (x\,\sigma_q-a)^{\alpha+1}\,\hat E_q^{(a,\alpha)}(x)=1\,,
\end{equation}
and is $Gq$-summable for any direction $d$ such that $a\notin(0,\infty e^{id})$. Furthermore, its associated $Gq$-sum is defined as follows:
\begin{equation}
 \label{equation:EaS}
 \mathcal{S}_{q;1}^d\hat E_q^{(a,\alpha)}(x)=\mathcal{L}_{q;1}^d\left(\frac{1}{(\xi-a)^{\alpha+1}}\right)(x)
\end{equation}
for all $x\in\tilde\CC^*$.

Moreover, the following identity holds in $\tilde\CC^*$: 
\begin{equation}\label{equation:SEadD}
 \mathcal{S}_{q;1}^d\hat E_q^{(a,\alpha)}(x)=\frac{1}{\alpha!}\,{\partial_ a^\alpha}\,\mathcal{S}_{q;1}^d\hat E_q^{(a,0)}(x)\,.
\end{equation}

\end{lem}

\begin{proof}
First, we shall prove \eqref{equation:EaEq} by induction on $\alpha$. Replacing $x$ with $\displaystyle\left(-\frac{x}{a}\right)$ in \eqref{equation:qEulerEquation} yields that 
$$\left(-\frac xa\,\sigma_q+1\right)\hat E_q(-\frac xa)=1.$$
With $\displaystyle\hat E_q^{(a,0)}(x)=-\frac{1}{a}\,\hat E_q(-\frac{x}{a})$ we find
\eqref{equation:EaEq} for $\alpha=0$. 
Furthermore, let $k\ge 1$, and suppose that  $(x\,\sigma_q-a)^k\,\hat E_q^{(a,k-1)}=1$, where $\alpha=k-1$. Expand this equation as follows:
$$
\sum_{j=0}^{k}\binom {k}j(-a)^{k-j}(x\,\sigma_q)^j\hat E_q^{(a,k-1)}(x)=1\,.
$$
By taking the differentiation with respect to $a$ in both sides in the above, using the identities $(k-j)\displaystyle\,\binom kj=\frac{(k-j)k!}{j! (k-j)!} =k\binom{k-1}j$ yields that
$$-
k\,\sum_{j=0}^{k-1}\binom {k-1}j(-a)^{k-j-1}\,(x\,\sigma_q)^j\hat E_q^{(a,k-1)}(x)+\sum_{j=0}^k\binom kj\,(-a)^{k-j}\,(x\,\sigma_q)^j{\partial_ a}\hat E_q^{(a,k-1)}(x)=0
$$
or, equivalently,
$$-k(x\,\sigma_q-a)^{k-1}\hat E_q^{(a,k-1)}(x)+(x\,\sigma_q-a)^{k}{\partial_ a}\hat E_q^{(a,k-1)}(x)=0\,.
$$
By \eqref{equation:qEulera}, one knows that $\displaystyle{\partial _a}\hat E_q^{(a,k-1)}=k\,\hat E_q^{(a,k)}$. One gets from the above that
$$(x\,\sigma_q-a)^{k}\,\hat E_q^{(a,k)}(x)=(x\,\sigma_q-a)^{k-1}\hat E_q^{(a,k-1)}(x)\,.
$$
Multiplying both sides by the operator $(x\,\sigma_q-a)$ and using $(x\,\sigma_q-a)^k\,\hat E_q^{(a,k-1)}=1$ implies \eqref{equation:EaEq} for $\alpha=k$.

By \cite[Equation (2.1)]{dreyfus2015building} with $\mu=1$, we find
\begin{equation}
 \label{equation:qBorelx}
 \hat\cB_{q;1}\,x^{j}\,\sq^{m}=q^{-j(j-1)/2}\xi^{j}\,\sq^{m-j}\hat\cB_{q;1}
\end{equation}
for all $(j,m)\in \ZZ_{\geq 0} \times\ZZ$. When  $j=m=1$, one gets that $\hat\cB_{q;1}(x\,\sigma_q)=\xi\,\hat\cB_{q;1}$. Thus, by applying the $q$-Borel transform $\hat\cB_{q;1}$ to both sides of \eqref{equation:EaEq}, it follows that $(\xi-a)^{\alpha+1}\hat\cB_{q;1}\hat E_q^{(a,\alpha)}=1$, {\it i.e.}
\begin{equation}
 \label{equation:EaBorel}
 \hat\cB_{q;1}\hat E_q^{(a,\alpha)}=\frac{1}{(\xi-a)^{\alpha+1}}\,.
\end{equation}
The right-hand side tends to $0$ as $|\xi|$ goes to $\infty$.  Hence, for any direction $d$ such that $a\notin(0,\infty e^{id})$,  $\displaystyle\frac{1}{(\xi-a)^{\alpha+1}}\in \HH_{q;1}^d$.
This implies the $Gq$-sum of $\hat E_q^{(a,\alpha)}$ as being stated in \eqref{equation:EaS}.

Given any couple of compact sets $K\subset\tilde\CC^*$ and $J\subset\CC^*$, the integral \eqref{equation:EaS} is uniformly convergent for $(x,a)\in K\times J$, provided that $J\cap(0,\infty\,e^{id})=\emptyset$. Thus, one can consider the differentiation with respect to $a$ under the integral sign. In view of the elementary identity 
$\displaystyle{\partial_ a^\alpha}\,\frac{1}{\xi-a}=\frac{\alpha!}{(\xi-a)^{\alpha+1}}$, one deduces easily \eqref{equation:SEadD} from \eqref{equation:EaS}. 
This finishes the proof of Lemma \ref{lemma:Ea}.
\end{proof}

By following \cite[Th\'eor\`eme 3.3 \& D\'efinition 2.1.1]{zhang1999developpements}, one knows that every $Gq$-sum $\mathcal{S}_{q;1}^d\hat E_q^{(a,\alpha)}$ admits the power series $\hat E_q^{(a,\alpha)}$ as a $q$-Gevrey asymptotic expansion of the first order in the direction $d$. It will be helpful to notice the following remark, due to the fact that
$
\left|e_q(t)\right|=e_q(|t|)\,e^{-(\arg(t))^2/(2\log (q))}$; see \eqref{equation:e_q} for the definition of $e_q$. 

\begin{rmk}\label{rmk:Eabouded}
 Let $a$, $\alpha$ and $d$ be as in Lemma \ref{lemma:Ea}. The $Gq$-sum  function $\mathcal{S}_{q;1}^d\hat E_q^{(a,\alpha)}$ is well-defined and analytic over the whole Riemann surface $\tilde\CC^*$, and satisfies the following inequality for all $x\in\tilde\CC^*$:
 \begin{equation}\label{equation:Eabounded}
 \left|\mathcal{S}_{q;1}^d\hat E_q^{(a,\alpha)}(x)\right|\le |a|_d^{-\alpha-1}\,e^{(\arg(e^{id}/x))^2/(2\log (q))}\,, 
 \end{equation}
where $|a|_d$ denotes the distance from $a$ to the half-line $(0,\infty e^{id})$ in the complex plane.
\end{rmk}

To see \eqref{equation:Eabounded}, one can notice that 
\begin{equation*}
\begin{split}
\left|\mathcal{S}_{q;1}^d\hat E_q^{(a,\alpha)}(x)\right| &\leq \int_0^{\infty}\frac{1}{|\xi e^{id}-a|^{\alpha+1}\, |e_q({\xi e^{id} }/{x})|}\,\frac{d\xi}{\xi}\\
&\leq \frac{e^{(\arg(e^{id}/x))^2/(2\log(q))}}{|a|_{d}^{\alpha+1}}\, \int_0^{\infty }\frac1{ e_q(\xi/|x|)}\,\frac{d\xi}{\xi}\,,
\end{split}
\end{equation*}
where the last integral is simply $\mathcal{L}_{q;1}^0 1(|x|)=1$.

 In what follows, we will make use of the following notation: 
 \begin{equation}\label{equation:Ea<}\forall k\ge0,\quad 
\hat E_q^{(a,-k-1)}=(-a)^k\,\sum_{j=0}^k\binom kj\,q^{j(j-1)/2}\,\left(-\frac xa\right)^j\,.
 \end{equation}
 Particularly, it follows that $\hat E_q^{(a,-1)}=1$.
 
Note that when $k<0$,  we do not recover $E_q^{(a,-k-1)}$. On the other hand, the notation fits to obtain the following lemma.

\begin{lem}\label{lemma:Eak}
 Let $(a,\alpha,k)\in\CC^*\times\ZZ_{\ge 0}\times\ZZ_{\ge0}$. Then:
 \begin{equation}\label{equation:Eak}
(x\,\sigma_q-a)^k\,\hat E_q^{(a,\alpha)}=\hat E_q^{(a,\alpha-k)}\,.
 \end{equation}

\end{lem}

\begin{proof}First, assume that $\alpha\ge k$.
 By applying \eqref{equation:EaEq}, both sides of the equation in \eqref{equation:Eak} are power series solutions of the linear $q$-difference equation $(x\,\sigma_q-a)^{\alpha-k+1}y=1$. Thus, their difference satisfies the homogeneous equation $(x\,\sigma_q-a)^{\alpha-k+1}y=0$, for which zero is the only power series solution. This implies \eqref{equation:Eak}.
 
 Next, when $k>\alpha$, we set $k'=k-\alpha-1 >0$, noticing that
 \begin{equation}\label{equation:Eak'}
 (x\,\sigma_q-a)^k\,\hat E_q^{(a,\alpha)}=(x\,\sigma_q-a)^{k'}\,\left((x\,\sigma_q-a)^{\alpha+1}\,\hat E_q^{(a,\alpha)}\right)=(x\,\sigma_q-a)^{k'}\,\cdot\,1\,.
 \end{equation}
 Since $(x\,\sigma_q)^j\,\cdot\,1=q^{j(j-1)/2}\,x^j\,\sigma_q^j\,\cdot\,1=q^{j(j-1)/2}\,x^j$ for all $j\in\ZZ$, the last term in \eqref{equation:Eak'} becomes:
 $$(x\,\sigma_q-a)^{k'}\,\cdot\,1=\sum_{j=0}^{k'}\binom{k'}j\,(-a)^{k'-j}\,q^{j(j-1)/2}\,x^j\,.
 $$Thus, 
 comparing \eqref{equation:Eak'} with  \eqref{equation:Ea<} yields \eqref{equation:Eak}.
 \end{proof}
 
 \subsection{Linear dependence of the generalized $q$-Euler series}\label{subsection:qEuler2}
 
 We shall start with the following elementary result on binomial coefficients.
 
 \begin{lem}\label{lemma:Ckj}
  If $\epsilon=0$ or $1$, then:
  \begin{equation}\label{equation:Ckj}
  \sum_{\ell=j}^{k+1-\epsilon}(-1)^{j+\ell}\,\binom {k+1}{\ell+\epsilon}\,\binom \ell j=\epsilon
  \end{equation}
  for any integers $k$ and $j$ such that $0\le j\le k$.
   \end{lem}
 
 \begin{proof}
  Let $C_{k,j}(\epsilon)$ be the left-hand side of \eqref{equation:Ckj}. 
  One has\begin{equation*}
  \begin{split}
  C_{k,j}(0)&=\sum_{\ell=j}^{k+1}(-1)^{j+\ell}\,\binom {k+1}{\ell}\,\binom \ell j=\binom{k+1}j\,\sum_{\ell=j}^{k+1}(-1)^{j+\ell}\,\binom {k+1-j}{\ell-j}\\
 &=\binom{k+1}j\,\sum_{r=0}^{k+1-j}(-1)^{r}\,\binom {k+1-j}{r}.
 \end{split}
  \end{equation*}
Since  $\displaystyle 0=(1-1)^{k+1-j}=\sum_{r=0}^{k+1-j}(-1)^{r}\,\binom {k+1-j}{r}$, one obtains that \eqref{equation:Ckj} for $\epsilon=0$. \par 
  
 Consider the case of $\epsilon=1$. It is easy to see that $C_{j,j}(1)=1$. Let $j \leq k$. Since
$$
 C_{k+1,j}(1)=\sum_{\ell=j}^{k}(-1)^{j+\ell}\,\binom {k+2}{\ell+1}\,\binom \ell j+(-1)^{j+k+1}\,\binom {k+2}{k+2}\binom {k+1}{j}\,,
 $$
  using the equalities $\displaystyle\binom{k+2}{\ell+1}=\binom{k+1}{\ell+1}+\binom{k+1}\ell$ and $\displaystyle\binom{k+2}{k+2}=\displaystyle\binom{k+1}{k+1}$ gives that
\begin{equation*}
  \begin{split}
 C_{k+1,j}(1)=& \displaystyle\sum_{\ell=j}^{k}(-1)^{j+\ell}\,\binom {k+1}{\ell+1}\,\binom \ell j\\&+\sum_{\ell=j}^{k}(-1)^{j+\ell}\,\binom {k+1}{\ell}\,\binom \ell j+(-1)^{j+k+1}\,\binom {k+1}{k+1}\,\binom {k+1}{j}\,.
 \end{split}
  \end{equation*}
  Thus, $
  C_{k+1,j}(1)=C_{k,j}(1)+C_{k,j}(0)=C_{k,j}(1)$ for  $0\le j\le k$.
In this way, one gets \eqref{equation:Ckj} for $\epsilon=1$. 
 \end{proof}

 \begin{thm}\label{theorem:Eaqk}
 Let $(a,\alpha,k)\in\CC^*\times\ZZ_{\ge 0}\times\ZZ_{\ge 0}$. One has
\begin{equation}\label{equation:Ea/qk}
\hat E_q^{(a/q^k,\alpha)}=q^{-k(k-2\alpha-3)/2}\,\left(\frac ax\right)^{k}\,\sum_{j=0}^k\binom kj\,\left(\frac1a\right)^j\,\hat E_q^{(a,\alpha-j)}
 \end{equation}
 and
  \begin{equation}\label{equation:Ea*qk}
   \begin{split}
\hat E_q^{(a,\alpha)}=
 q^{k(k-2\alpha-3)/2} \left(\frac x{ a}\right)^k\,\sum_{\ell=0}^\alpha\binom{k+\ell-1}{k-1}\,\left(-\frac{q^k}a\right)^\ell\,\hat E_q^{(a/q^k,\alpha-\ell)}\\
 +\left(-\frac{1}a\right)^{\alpha+1}\,\sum_{j=0}^{k-1}q^{j(j-1)/2}\,\binom{\alpha+j}{\alpha}\,\left(\frac xa\right)^{j}\,.
  \end{split}
 \end{equation}
\end{thm}

\begin{proof} 
 Replacing $\lambda$ with $q^k$ in \eqref{equation:qEuleraa'} yields that $\hat E_q^{(a/q^k,\alpha)}=q^{(\alpha+1)k}\,\sigma_q^{k}\,\hat E_q^{(a,\alpha)}$.
 By multiplying this last equation by $q^{k(k-1)/2}\,x^k$ and observing that $(x\,\sigma_q)^k=q^{k(k-1)/2}\,x^k\,\sigma_q^k$, one finds that
 \begin{equation}\label{equation:Eaqk}
(x\,\sigma_q)^{k}\,\hat E_q^{(a,\alpha)}=q^{k(k-2\alpha-3)/2}\,x^k\,\hat E_q^{(a/q^k,\alpha)}\,.
 \end{equation}
 On the other hand, by using the identity
 $$
 (x\,\sigma_q)^{k}=\left((x\,\sigma_q-a)+a\right))^{k}=\sum_{j=0}^k\binom kj\,a^{k-j}\,(x\,\sigma_q-a)^{j}\,,
 $$
 one deduces from gathering \eqref{equation:Eaqk} together with \eqref{equation:Eak} the identity stated in \eqref{equation:Ea/qk}.
 
 Now, let $\alpha=0$ into \eqref{equation:Ea/qk}, and separe $j=0$ and $j=\ell+1>0$ in the summation. This implies that 
 \begin{equation}\label{equation:EaqkA}
\hat E_q^{(a/q^k,0)}=q^{-k(k-3)/2}\,\left(\frac ax\right)^{k}\,\left(\hat E_q^{(a,0)}+A_k\right)\,,
 \end{equation}
 where $A_0=0$ and for $k\ge 1$:
$$
A_k=\sum_{\ell=0}^{k-1}\binom k{\ell+1}\,\left(\frac1a\right)^{\ell+1}\,\hat E_q^{(a,-\ell-1)}\,.
$$Assume $k>0$, and
apply
\eqref{equation:Ea<} for each term $\hat E_q^{(a,-\ell-1)}$ in the above. One gets that
$$
A_k=\frac1a\,\sum_{\ell=0}^{k-1}\sum_{j=0}^\ell(-1)^\ell\,\binom k{\ell+1}\,\binom \ell j\,q^{j(j-1)/2}\,\left(-\frac xa\right)^j\,,
$$
or, equivalently,
$$
A_k=\frac1a\,\sum_{j=0}^{k-1}\sum_{\ell=j}^{k-1}(-1)^{j+\ell}\,\binom k{\ell+1}\,\binom \ell j\,q^{j(j-1)/2}\,\left(\frac xa\right)^j\,.
$$
By considering \eqref{equation:Ckj} for $\epsilon=1$ and $(k-1)$ instead of $k$, it follows that
$$ A_k
=\frac1a\,\sum_{j=0}^{k-1}q^{j(j-1)/2}\,\left(\frac xa\right)^j\,.
$$
 So, \eqref{equation:EaqkA} can be put into the following form:
 \begin{equation}\label{equation:qEulerak}
  \hat E_q^{(a,0)}(x)=q^{k(k-3)/2}\left(\frac xa\right)^k\,\hat E_q^{(a/q^k,0)}(x)-\frac{1}{a}\,\sum_{j=0}^{k-1}q^{j(j-1)/2}\,\left(\frac{x}{a}\right)^j\,,
 \end{equation}
 which is really \eqref{equation:Ea*qk} for $\alpha=0$.

 To obtain \eqref{equation:Ea*qk} for $\alpha>0$, we will apply the definition of $\hat E_q^{(a,\alpha)}$ in \eqref{equation:qEulera}, considering the $\alpha$-th differential of $\hat E_q^{(a,0)}$ with the help of \eqref{equation:qEulerak}. Indeed, given a couple of non-negative integers $(\ell,j)$, one has
 $$\partial_a^\ell\left(a^{-j-1}\right)=(-1)^\ell\,(j+1)\,(j+2)\,\dots\,(j+\ell)\,a^{-j-\ell-1}
 $$
 and, in addition, using \eqref{equation:qEulera} implies that 
 $
 \partial_a^\ell\hat E_q^{(a/q^j,0)}=\ell!\,\hat E_q^{(a/q^j,\ell)}\,q^{-j\,\ell}$. 
 Thus, by taking the $\alpha$-th differential with respect to $a$ for both sides of \eqref{equation:qEulerak}, it follows  that
 \begin{equation*}
  \begin{split}
\alpha!\,\hat E_q^{(a,\alpha)}=
 q^{k(k-3)/2} \left(\frac xa\right)^k\,\sum_{\ell=0}^\alpha\binom\alpha\ell\,\frac{(k+\ell-1)!}{(k-1)!}\,\left(-\frac1a\right)^\ell\,(\alpha-\ell)!\,q^{-k(\alpha-\ell)}\,\hat E_q^{(a/q^k,\alpha-\ell)}\\
 -\sum_{j=0}^{k-1}q^{j(j-1)/2}\,x^j\,(-1)^\alpha\,\frac{(j+\alpha)!}{j!}\,a^{-j-\alpha-1}\,.
  \end{split}
 \end{equation*}
 By dividing by $\alpha!$ in the above, one gets the expression given in \eqref{equation:Ea*qk}.
\end{proof}

\begin{cor}\label{corollary:EaCL}
Let $a\in\CC^*$, $m\in\ZZ_{>0}$, $(k_1,\dots,k_m)\in\ZZ_{\ge0}^m$,  $(\alpha_1,\dots,\alpha_m)\in\ZZ_{\ge 0}^m$, and let $(f_1, \dots, f_m)\in\left(\CC\{x\}\right)^m$. If $k=\max(k_1,\dots,k_m)$, $a'=a/q^k$ and $\alpha=\max(\alpha_1,\dots ,\alpha_m)$, then one can find $(H,h_0,h_1, \dots , h_\alpha)\in\left(\CC\{x\}\right)^{\alpha+2}$ such that
\begin{equation}\label{equation:EaCL}
 \sum_{j=1}^mf_j\,\hat E_q^{(a/q^{kj},\alpha_j)}=H+\sum_{j=0}^\alpha h_j\,\hat E_q^{(a',\alpha-j)}\,.
\end{equation}

\end{cor}

\begin{proof}
 This comes directly from applying \eqref{equation:Ea*qk}.
\end{proof}

 \subsection{Generic irregular singular $q$-difference operators coming from combining several generalized $q$-Euler series}\label{subsection:GenericOperators}

 For simplicity, let $\DD_0=\CC^*\oplus\left(x\,\CC\{x\}\right)$, let $\DD _n$ be the set of the $n$-th order $q$-operators $\Delta$ given in \eqref{equation:Delta} for $n\ge 1$, and define
 \begin{equation}\label{equation:DD}
\DD =\cup_{n\ge 0}\DD _n\,.
 \end{equation}
 We call {\it a generic irrgular singular $q$-difference operator} any element $\Delta\in\DD_n$ with $n\ge 1$.
 
 By following \cite[Prop. 5.1.4]{zhang1999developpements}, each operator $\Delta\in\DD_n $ can be factorized as follows:
 \begin{equation}\label{equation:fact}
  \Delta=h_0\,(x\,\sigma_q-a_1)\,h_1\,\dots \,(x\,\sigma_q-a_n)\,h_n\,,
 \end{equation}
 where $a_j\in\CC^*$ and $h_j\in\DD_0$. Let  $c\mapsto\bar c$ be the natural map from $\CC^*$ to the multiplcative elliptic curve $\CC^*/q^\ZZ$. Define the set of (generalized) indices of monodromy $I\!M(\Delta)$ as follows: $I\!M(h)=\emptyset$ for $h\in\DD_0$, and
 \begin{equation}\label{equation:IM}
I\!M(\Delta)=\left\{\bar a_1,\bar a_{2},\dots,\bar a_n\right\}\subset\CC^*/q^\ZZ
 \end{equation}
 for $\Delta\in\DD_n$, $n\ge 1$.
  One knows that $I\!M(\Delta)$ consitutes a set of analytic invariants of $\Delta$; for instance, see \cite[\S 5.1.7]{zhang1999developpements}.
In addition, given $(\Delta_1,\Delta_2)\in \DD _{n_1}\times \DD _{n_2}$, one has
 \begin{equation}\label{equation:DD1}
 \Delta_1\,\Delta_2\in\DD _{n_1+n_2}\,,\quad I\!M(\Delta_1\,\Delta_2)=I\!M(\Delta_1)\cup I\!M(\Delta_2)\,.
 \end{equation}
 
Furthermore, for any pair  $(\Delta,h)\in\DD_n\times\left(\CC\{x\}\setminus\{0\}\right)$ ($n\ge 1$), if $\Delta^{[h]}$ denotes the conjugation of $\Delta$ by $h$, {\it i.e.}
$\Delta^{[h]}=(\sigma_q^nh)^{-1}\,\Delta\,h$, where $(\sigma_q^nh)^{-1}(x)=\displaystyle\frac{1}{h(q^n\,x)}$, one has the following properties: 
\begin{equation}\label{equation:Dh}
\Delta^{[h]}\in\DD _n\,,\quad 
 I\!M(\Delta^{[h]})=I\!M(\Delta)
\end{equation}
and \begin{equation}\label{equation:Dhy}
\Delta\,(h\,y)=\sigma_q^nh\,\Delta^{[h]}y\,.
\end{equation}

\begin{lem}\label{lemma:Euclid}
 Given $(a,n,\Delta)\in\CC\times\ZZ_{>0}\times\DD_n$,  there exist $b_0, b_1, \dots , b_n\in\CC\{x\}$ such that the following expression holds:
 \begin{equation}\label{equation:Euclid}
 \Delta=b_0\,(x\,\sigma_q-a)^n+b_1\,(x\,\sigma_q-a)^{n-1}+\dots+b_n\,.  
 \end{equation}

\end{lem}

\begin{proof} 
Let $a_0, \dots, a_n$ be the coefficients of $\Delta$ as given in \eqref{equation:Delta}. By writing $x\,\sigma_q=(x\,\sigma_q-a)+a$, applying the binomial expansion for $(x\,\sigma_q)^ k$ yields that
 $$(x\,\sigma_q)^k=\sum_{j=0}^k\binom kj\,a^{k-j}\,(x\,\sigma_q-a)^{j}\,.$$
Putting this into \eqref{equation:Delta} yields that 
$$\Delta=\sum_{k=0}^{n}a_{n-k}\sum_{j=0}^k\binom kj\,a^{k-j}\,(x\,\sigma_q-a)^{j}\,.$$
Therefore, if
\begin{equation}\label{equation:Euclidb}
 b_j=
\sum_{\ell=0}^j\binom {n-\ell}{n-j}\,a^{j-\ell}\,a_{\ell}\quad(0\le j\le n)\,,
 \end{equation}
 one gets immediately \eqref{equation:Euclid} and finishes the proof of Lemma \ref{lemma:Euclid}.
\end{proof}

\begin{rmk}\label{remark:Euclid}
The above relations in
\eqref{equation:Euclidb} lead us to the following elementary identity:
\begin{equation}\label{equation:Euclida}
\forall a\in\CC\,,\quad \CC\{x\}[x\,\sigma_q]=\CC\{x\}[x\,\sigma_q-a]\,.
\end{equation}
\end{rmk}

\begin{thm}\label{theorem:fEaEq}
Let $(a,\alpha, m)\in \CC^*\times\ZZ_{\ge 0}\times\ZZ_{>0}$, and let $(\hat f, \hat g,h_0,h_1,\dots,h_{\alpha}, \Lambda )\in\CC[[x]]\times\CC[[x]]\times\left(\CC\{x\}\right)^{\alpha+1}\times\DD_{m }$ to be  such that 
\begin{equation}\label{equation:fgEa}
\hat f=\hat g+\sum_{\ell=0}^\alpha h_{\alpha-\ell}\,\hat E_q^{(a,\ell)}\,,\quad \Lambda\, \hat g\in\CC\{x\}\,.
\end{equation}
One supposes that $\bar a\notin I\!M(\Lambda )$ and $h_0\not=0$. Then, there exists $\Delta\in\DD_{n}$ such that $0\leq n\le \frac12\left({\alpha+1}\right)\left({\alpha+2}\right)$,
\begin{equation}\label{equation:fEaEq}
\Delta\,\Lambda \hat f\in\CC\{x\},\,\quad I\!M(\Delta)=\{\bar a\}\,.
\end{equation}

\end{thm}

\begin{proof} 
By considering \eqref{equation:Dhy}, one knows that 
$$\Lambda (h_{\alpha-\ell}\,\hat E_q^{(a,\ell)})=\sigma_q^{m }h_{\alpha-\ell}\,\Lambda ^{[h_{\alpha-\ell}]}\hat E_q^{(a,\ell)}$$
for $\ell$ from 0 to $\alpha$.
 Applying Lemma \ref{lemma:Euclid} to $\Lambda ^{[h_{\alpha-\ell}]}$ gives that
\begin{equation}\label{equation:Delta1}
 \Lambda ^{[h_{\alpha-\ell}]}=\sum_{j=0}^{m }b_{\alpha-\ell,j}\,(x\,\sigma_q-a)^{m -j}\,,\quad b_{\alpha-\ell,j}\in\CC\{x\}\,.
\end{equation}
From now on, define
\begin{equation}\label{equation:Yj}h_{\ell,j}=b_{\alpha-\ell,m-j}\,\sigma_q^{m }h_{\alpha-\ell}\,,\quad
 \hat Y_{\ell,j}=(x\sigma_q-a)^{j}\,\hat E_q^{(a,\ell)}\,.
\end{equation}
In view of \eqref{equation:EaEq}, it follows that $\hat Y_{\ell,j}\in\CC\{x\}$ when $j>\ell$. Moreover, using \eqref{equation:Eak} gives that $\hat Y_{\ell,j}=\hat Y_{\ell-j,0}=\hat E_q^{(a,\ell-j)}$ if $j\le\ell$.
So, applying $\Lambda $ to $\hat f$ in \eqref{equation:fgEa} and gathering \eqref{equation:Delta1} together with \eqref{equation:Yj} implies that
$$
\Lambda \hat f=H_0+\sum_{\ell=0}^\alpha\sum_{j=0}^{\ell }h_{\ell,j}\,\hat E_q^{(a,\ell-j)}\,,\quad H_0\in\CC\{x\}\,.
$$
Letting $\displaystyle h_k^*=\sum_{j=0}^{\alpha-k }h_{j+k,j}$, it follows from the above that
\begin{equation}\label{equation:Yj1}
\Lambda \hat f=H_0+\sum_{k=0}^\alpha h_k^*\,\hat E_q^{(a,k)}\,,\quad H_0\,,\ h_k^*\in\CC\{x\}\,.
\end{equation}

In what follows, we will construct a family of $q$-difference operators to successively remove from  \eqref{equation:Yj1} the divergent power series $\hat E_q^{(a,0)}$, $\hat E_q^{(a,1)}, \dots, \hat E_q^{(a,\alpha)}$. For doing this, define 
\begin{equation}\label{equation:hDelta}
\Delta_{0}=
\left\{ 
\begin{array}{ll}
\left((x\,\sigma_q-a)^{\alpha+1}\right)^{[1/h_{\alpha}^*]}&\hbox{if $h_{\alpha}^*\not=0$};\cr
1&\hbox{if $h_{\alpha}^*=0$}.\cr
\end{array}
\right. 
\end{equation}
By applying \eqref{equation:Dhy} to $\displaystyle\left(\Delta,h\right)=\left((x\,\sigma_q-a)^{\alpha+1},1/h_{\alpha}^*\right)$, it follows from considering \eqref{equation:EaEq}  that, in any case where $h_{\alpha}^*=0$ or not,
$$\Delta_{0}(h_\alpha^*\,\hat E_q^{(a,\alpha)})=\sigma_q^{\alpha+1}\,h_{\alpha}^*\,(x\,\sigma_q-a)^{\alpha+1}\,\hat E_q^{(a,\alpha)}=\sigma_q^{\alpha+1}\,h_\alpha^*\,.
$$
In this way, 
one obtains from \eqref{equation:Yj1} that
\begin{equation}\label{equation:hDelta2}
\Delta_{0}\,\Lambda \hat f=H_1+\sum_{k=0}^{\alpha-1}\Delta_0\left(h_{k}^*\,\hat E_q^{(a,k)}\right)\,,
\end{equation}
where $H_1=\Delta_{0}\,H_0+\sigma_q^{\alpha+1}h_\alpha^*$. We remark that the sum is zero when $\alpha=0$. Consequently, letting $\Delta=\Delta_0$ gives \eqref{equation:fEaEq} if $\alpha=0$.

If $\alpha>0$, again by applying \eqref{equation:Dhy}, one has $$\Delta_{0}\,(h_{k}^*\,\hat E_q^{(a,k)})=\sigma_q^{\alpha+1}h_{k}^*\,\Delta_{0}^{[h_{k}^*]}\,\hat  E_q^{(a,k)}$$ for $0\le k< \alpha$. One deduces from \eqref{equation:hDelta2} that
\begin{equation}\label{equation:hDelta1}
\Delta_{0}\,\Lambda \hat f=H_1+\sum_{k=0}^{\alpha-1}\sigma_q^{\alpha+1}h_{k}^*\,\Delta_{0}^{[h_{k}^*]}\,\hat E_q^{(a,k)}\,.
\end{equation}
By expanding each operator $\Delta_{0}^{[h_{k}^*]}$ in $\CC\{x\}[x\,\sigma_q-a]$ (see Lemma \ref{lemma:Euclid}) and observing that $(x\,\sigma_q-a)^j\,\hat E_q^{(a,k)}=\hat E_q^{(a,k-j)}$  for $k\ge j$ (see \eqref{equation:Eak}), one can put \eqref{equation:hDelta1} into the following form:
\begin{equation}\label{equation:hDelta1*}
\Delta_{0}\,\Lambda \hat f=H_1^*+\sum_{k=0}^{\alpha-1} h_{k}^{**}\,\hat E_q^{(a,k)}\,,\quad H_1^*\,,\ h_k^{**}\in\CC\{x\}\,.
\end{equation}

It is to be noted that the divergent series $\hat E_q^{(a,\alpha)}$ does not appear in the right-hand side of the equation in \eqref{equation:hDelta1*}. By a similar process as what done for the transformation from \eqref{equation:Yj1} into \eqref{equation:hDelta1*}, one can define a series of operators $\Delta_{1}$, $\Delta_{2}, \dots, \Delta_{\alpha}$ possessing the following properties:
\begin{enumerate}
 \item 
either $\Delta_{k}\in\DD_{k+1}$ and  $I\!M(\Delta_{k})=\{\bar a\}$ or $\Delta_{k}=1$;
\item letting $\Delta=\Delta_{\alpha}\,\dots\,\Delta_{1}\,\Delta_{0}$, one has $\Delta\,\Lambda \hat f\in\CC\{x\}$. 
\end{enumerate}
Since the order of $\Delta$ equals at most to the following sum: 
$$
\sum_{k=0}^\alpha(\alpha+1-k)=\frac12(\alpha+1)(\alpha+2)\,,
$$
one gets \eqref{equation:fEaEq} and finishes the proof of Theorem \ref{theorem:fEaEq}. 
\end{proof}

\subsection{Characterisation of generic $q$-Gevrey series and proof of Theorem \ref{theorem:generic}}\label{subsection:thmGeneric}
As in \cite{di2009q}, let $K$ be the set of rational functions bounded at zero, and let $\EE_{q;1}$ be the set of all entire functions admitting at most a $q$-exponential growth of order 1 at the infinity. It follows that 
\begin{equation}\label{equation:4E}
\EE_{q;1}=\cap_{d\in[0,2\pi]}\HH_{q;1}^d=\hat\cB_{q;1}(\CC\{x\})\,. 
\end{equation}

Let $\star_{q;1}$ be the bilinear map from $\CC[[\xi]]\times \CC[[\xi]]$ to $\CC[[\xi]]$ such that  $\xi^m\star_{q;1}\xi^n=q^{-m\,n}\,\xi^{m+n}$ for all $(m,n)\in\ZZ_{\ge 0}^2$. Thus, the $\CC$-vector space $\CC[[\xi]]$ is equiped with two {\it multiplications} that are the usual one and the so-called {\it $q$-convolution product} $\star_{q;1}$. Furthermore,  the following identity 
holds for any $(\hat f, \hat g)\in\CC[[x]]\times\CC[[x]]$:
\begin{equation}\label{equation:4star}
\hat\cB_{q;1}(\hat f\,\hat g)=\hat\cB_{q;1}\hat f\star_{q;1}\hat\cB_{q;1}\hat g\,.
\end{equation}

By following \cite[\S 4]{di2009q}, one defines the $(K,\EE_{q;1})$-bimodule $H\subset\CC\{\xi\}\subset\CC[[\xi]]$ by the relation $H=\cup_{n\ge 0}H_n$, where $ H_0=K$
and, for any integer $n\ge 0$,
$$
H_{2n+1}=\EE_{q;1}\star_{q;1} H_{2n}\,,\quad H_{2n+2}=K\,H_{2n+1}\,.
$$ 
By \cite[Theorem 4.20]{di2009q}, one knows that $\phi\in H$ if and only if 
\begin{equation}\label{equation:4phi}
\phi=r_0+\phi_0+\phi_1\star_{q;1}\frac{1}{(\xi-a_1)^{\alpha_1+1}}+\dots+\phi_n\star_{q;1}\frac{1}{(\xi-a_n)^{\alpha_n+1}}\,,
\end{equation}
where $r_0\in K$, $\phi_j\in\EE_{q;1}$, $a_j\in\CC^*$ and $\alpha_j\in\ZZ_{\ge 0}$. 
By noticing that both $K$ and $\EE_{q;1}$ are stable by $\sigma_q$ and that, in addition, either the multiplication or the $q$-convolution product is compatible with $\sigma_q$, one finds that $H$ is a  $q$-difference $(K,\EE_{q;1})$-bimodule such that
\begin{equation}\label{equation:4phiq}
 \phi\in H\quad\Longleftrightarrow\quad\sigma_q\,\phi\in H\,.
\end{equation}

\begin{lem}\label{lemma:4phi}
 Given any $\hat f\in\CC[[x]]$, one has $\hat\cB_{q;1}\hat f\in H$ if and only if $\hat f $ can be put into the form of \eqref{equation:decomposition}.
\end{lem}
\begin{proof}
 First of all, suppose $\hat f$ be as in  \eqref{equation:decomposition}. By considering \eqref{equation:4star}, one finds that
 \begin{equation}\label{equation:4phif}
  \hat\cB_{q;1}\hat f=\hat\cB_{q;1} f_0+\hat\cB_{q;1} f_1\star_{q;1}\hat\cB_{q;1}\hat E_q^{(a_1,\alpha_1)}+\dots+\hat\cB_{q;1} f_n\star_{q;1}\hat\cB_{q;1}\hat E_q^{(a_n,\alpha_n)}\,.
 \end{equation}
 Furthermore, by applying \eqref{equation:EaBorel}, one has $\displaystyle \hat\cB_{q;1}\hat E_q^{(a_j,\alpha_j)}=\frac{1}{(\xi-a_j)^{\alpha_j+1}}$. Thus, letting $\phi=\hat \cB_{q;1}\hat f$, $r_0=0$ and  $\phi_j=\hat\cB_{q;1} f_j$ for $0\le j\le n$, one deduces \eqref{equation:4phi} from \eqref{equation:4phif}. This implies that $\hat\cB_{q;1}\hat f\in H$.
 
Next, suppose $\hat\cB_{q;1}\hat f\in H$, and assume \eqref{equation:4phi} hold for $\phi=\hat\cB_{q;1}\hat f$. By making use of the fractional decomposition of $r_0$, one can suppose that
$$
r_0=P+\sum_{\ell=1}^m\sum_{k=0}^{\mu_\ell-1}\frac{c_{\ell,k}}{(\xi-b_\ell)^{k+1}}\,,\quad P\in\CC[\xi]\,,\ c_{\ell,k}\in\CC,\ b_\ell\in\CC^*\,.
$$
Therefore, by following \eqref{equation:EaBorel}, there exists $Q\in\CC[x] $ such that  $$
r_0=\hat \cB_{q;1}\left(Q+\sum_{\ell=1}^m\sum_{k=0}^{\mu_\ell-1}c_{\ell,k}\,\hat E_q^{(b_\ell,k)}\right)\,.
$$
Thus, by a similar way as what done in the above,   one gets \eqref{equation:decomposition} from considering  \eqref{equation:4phi}, using \eqref{equation:4star} at the same time. We omit the detail.
\end{proof}

\begin{proof}[Proof of Theorem \ref{theorem:generic}] Let $\hat f=c_0+c_1\,x+c_2\,x^2+\dots\in\CC[[x]]$. We will proceed in two steps.

On the one hand, suppose $\hat f$ be a generic $q$-Gevrey series, that is to say, there exists $\Delta\in\DD$ such that $\Delta\hat f\in\CC\{x\}$. Set $\phi=\hat\cB_{q;1}\hat f$, $\hat g=\hat f-c_0$, and observe that $\Delta\hat g=\Delta\hat f-\Delta c_0\in\CC\{x\}$.
Now, consider the $q$-Borel transform  $B_q$ defined in \cite[p. 382]{di2009q}, which is the linear application $x\,\CC[[x]]\to\CC[[\xi]]$ such that $B_q(x^{n+1})=q^{-n(n-1)/2}\,\xi^n$ for all $n\in\ZZ_{\ge0}$. If $\gamma=B_q\hat g$, one finds that 
$\displaystyle
\gamma=\hat\cB_{q;1}\left(({\hat f-c_0})/x\right)$. By taking into account the relation $\hat f=c_0+x\,\left(({\hat f-c_0})/x\right)$, applying $\hat\cB_{q;1}$ to $\hat f$ and using \eqref{equation:4star} gives that
$$
\phi=c_0+\hat\cB_{q;1}x\star_{q;1}\hat\cB_{q;1}\left(\frac{\hat f-c_0}x\right)=c_0+\xi\star_{q;1}\gamma\,.
$$
Thus, the $q$-Borel transforms $\phi$ and $\gamma$ are linked as follows:
\begin{equation}\label{equation:4phi_gamma}
\phi=c_0+\xi\,\sigma_q^{-1}\gamma\,.
\end{equation}

By applying \cite[Proposition 4.24]{di2009q}, one knows that $\gamma\in H$. This together with considering both \eqref{equation:4phiq} and \eqref{equation:4phi_gamma} implies that $\phi\in H$. By Lemma \ref{lemma:4phi}, it follows that $\hat f$ can be put into the form as in \eqref{equation:decomposition}.

On the other hand,  suppose $\hat f$ be such that the relation in \eqref{equation:decomposition} holds, {\it i.e.}
$$
\hat f=f_0+f_1\,\hat E_q^{(a_1,\alpha_1)}+\dots+f_m\,\hat E_q^{(a_m,\alpha_m)}\,.
$$
Let $\mu$ be the number of different elements contained in the subset $\{\bar a_1,\dots,\bar a_m\}$ of the elliptic curve $\CC^*/q^\ZZ$. If $\mu=0$, then $\hat f=f_0\in\CC\{x\}$. As the convergent power series $\hat f$ is a generic $q$-Gevrey series, one can finish the proof of Theorem \ref{theorem:generic} by induction on $\mu$, using Theorem \ref{theorem:fEaEq}.
 \end{proof}

\section{The product of two generalized $q$-Euler series}\label{section:qEuler}

 This section aims to prove one special case of Theorem \ref{theorem:intro} for two generalized $q$-Euler  series. Namely, we shall establish the following statement.

\begin{thm}\label{theorem:Eab}
 Given $(a,b)\in\CC^*\times\CC^*$, one has $\hat E_q^{(a,0)}\,\hat E_q^{(b,0)}\in\CC\{x\}_{q;(1,2)}$. 
 Moreover, for any $d\in\RR$ such that $\{a,b\}\cap(0,\infty e^{id})=\emptyset$, one has
 \begin{equation}
  \label{equation:EaEbS}
  \mathcal{S}_{q;(1,2)}^d(\hat E_q^{(a,0)}\,\hat E_q^{(b,0)})=\mathcal{S}_{q;1}^d(\hat E_q^{(a,0)})\,\mathcal{S}_{q;1}^d(\hat E_q^{(b,0)})\,.
 \end{equation}

\end{thm}

The proof of Theorem \ref{theorem:Eab} will be given in \S \ref{subsection:proofThmEab}. To do that, we will start with a remarkable relation for special values of the $Gq$-sum functions of two generalized $q$-Euler series. This will permit us to obtain the $Gq$-summability of their product, after having considered their $q$-Borel transform and a first $q$-Laplace transform. 

\subsection{Special values of a combination of two generalized $q$-Euler series}\label{subsection:qEulerValues}
By following \eqref{equation:EaS}, one knows that, given $d\in\RR$ and $\alpha\in\ZZ_{\ge 0}$, the associated sum-function $\cS_{q;1}^d(\hat E_{q}^{(a,\alpha)})$ represents an analytic function with respect to the parameter $a\in\CC^*$, provided that $a\notin(0,e^{id})$. In Lemma \ref{lemma:EaEb} below, we shall consider sum-functions at some special points related to the parameter $a$. To simplify, we will make use of  the following conventional hypothesis.

\begin{hypo}
\label{hypo:abd}
The triplet $(a,b,d)\in\CC^*\times\CC^*\times\RR$ satisfies the following conditions:
 \begin{equation}
  \label{equation:abd}
 \left\{\arg(a),\arg(b)\right\}\subset(d,d+2\,\pi)\,.
 \end{equation}

\end{hypo}

Note that \eqref{equation:abd} above  implies 
 $\{a,b\}\cap (0,\infty e^{id})=\emptyset$.

\begin{lem}\label{lemma:EaEb}
Let $(a,b,d)\in\CC^*\times\CC^*\times\RR$ to be such that $(0,\infty e^{id})\cap\{a,b\}=\emptyset$ and that \eqref{equation:abd} holds,  and let $e_q$ be as in  \eqref{equation:e_q}. 
 Then, the following equality holds for any $\ell\in\ZZ$:
 \begin{equation}\label{equation:EaEb}
  a\,\mathcal{S}_{q;1}^d\hat E_q^{(a,0)}(\sqrt{q\,a\,b}\,e^{\pi i\ell})+b\,\mathcal{S}_{q;1}^d\hat E_q^{(b,0)}(\sqrt{q\,a\,b}\,e^{\pi i\ell})=-1+2\pi i\,R_{q;\ell}^d(a,b)\,,
 \end{equation}
where
\begin{equation}
 \label{equation:EaEbR}
 R_{q;\ell}^d(a,b)=\left\{
 \begin{array}{ll}\displaystyle-\sum_{j=\ell+1}^{-1}\frac{1}{e_q(\sqrt{b/(q\,a)}\,e^{(\ell-2j)\pi i})}&\hbox{for}\ \ell<-1\,;\cr
 0&\hbox{for}\ \ell=-1\,;\cr
  \displaystyle\sum_{j=1}^{\ell+1}\frac{1}{e_q(\sqrt{b/(q\,a)}\,e^{(\ell-2j)\pi i})}&\hbox{for}\ \ell>-1\,.\cr
 \end{array}
 \right.
\end{equation}

\end{lem}
\begin{proof} 
Let $c=\sqrt{q\,a\,b}$, $c_\ell=c\,e^{\pi i\ell}$, $a=|a|\,e^{i{a'}}$, $b=|b|\,e^{i{b'}}$, and let denote by $A_\ell$ the left-hand side of \eqref{equation:EaEb}.
By raplacing $(a,\alpha,x)$ with $(a,0,c_\ell)$ and $(b,0,c_\ell)$ in \eqref{equation:EaS} respectively, it follows from using \eqref{equation:qLaplace1} that
\begin{equation}
 \label{equation:EaEbA}A_\ell=\mathcal{L}_{q;1}^d\left(\frac{a}{\xi-a}+\frac{b}{\xi-b}\right)(c_\ell)=I_{1,\ell}+I_{2,\ell},
\end{equation}
where
$$
I_{1,\ell}=\int_0^{\infty e^{id}}\frac{1}{(\xi/a)-1}\,\frac{1}{e_q(\xi/c_\ell)}\,\frac{d\xi}{\xi}\,,\quad I_{2,\ell}=\int_0^{\infty e^{id}}\frac{1}{(\xi/b)-1}\,\frac{1}{e_q(\xi/c_\ell)}\,\frac{d\xi}{\xi}\,.
$$
Letting $t=\displaystyle\frac\xi a$ and $t=\displaystyle\frac\xi b$  respectively in $I_{1,\ell}$ and $I_{2,\ell}$ in the above gives that 
\begin{equation}
 \label{equation:EaEbI1}
I_{1,\ell}=\int_0^{\infty e^{i(d-{a'})}}\frac{1}{t-1}\,\frac{1}{e_q(a\,t/c_\ell)}\,\frac{dt}{t}
\end{equation}
and 
\begin{equation}
 \label{equation:EaEbI2} I_{2,\ell}=\int_0^{\infty e^{i(d-{b'})}}\frac{1}{t-1}\,\frac{1}{e_q(b\,t/c_\ell)}\,\frac{dt}{t}\,.
\end{equation}
By considering
\eqref{equation:e_q}, one knows that $e_q(x)=e_q(1/(qx))$. Thus, the integral in \eqref{equation:EaEbI1} can be rewritten as follows:
$$
I_{1,\ell}=\int_0^{\infty e^{i(d-{a'})}}\frac{1}{t-1}\,\frac{1}{e_q(c_\ell/(q\,a\,t))}\,\frac{dt}{t}
=\int_0^{\infty e^{i(a'-d)}}\frac{t}{1-t}\,\frac{1}{e_q(c_\ell\,t/(q\,a))}\,\frac{dt}{t}\,,
$$
the last integral being obtained by changing $t$ with $1/t$. By using the fractional decomposition $\displaystyle\frac{t}{1-t}=\frac{1}{1-t}-1$ and the relation   $\mathcal{L}_{q;1}^d(1)=1$, one finds that 
\begin{equation}
 \label{equation:EaEbI1a}I_{1,\ell}=-1+\int_0^{\infty e^{i(a'-d)}}\frac{1}{1-t}\,\frac{1}{e_q(c_\ell\,t/(q\,a))}\,\frac{dt}{t}\,.
\end{equation}

Consider the case $\ell=-1$, and set $c'=c_{-1}$, {\it i.e.} $c'=c\,e^{-\pi i}$. By taking into account the hypothesis made on the triplet $({a'},{b'},d)$, one finds that $\displaystyle\frac{b}{c'}=\frac{c'}{q\,a}\,e^{2\pi i}$. Thus, replacing $(\ell,t)$ with $(-1,t\,e^{-2\pi i})$ in \eqref{equation:EaEbI2}  yields that
$$
 I_{2,-1}=\int_0^{\infty e^{i(d-{b'}+2\pi)}}\frac{1}{t-1}\,\frac{1}{e_q(c'\,t/(q\,a))}\,\frac{dt}{t}\,.
$$
Gathering this  together with \eqref{equation:EaEbI1a} and \eqref{equation:EaEbA} yields that
\begin{equation}
 \label{equation:EaEbA0}
A_{-1}=-1+\left(\int_0^{\infty e^{i(a'-d)}}-\int_0^{\infty e^{i(d-{b'}+2\pi)}}\right)\frac{1}{1-t}\,\frac{1}{e_q(c'\,t/(q\,a))}\,\frac{dt}{t}\,.
\end{equation}
Furthermore, in view of the relation $\{{a'}-d,d-{b'}+2\pi\}\subset (0,2\,\pi)$, the integration-contour composed of the half-lines $(0,\infty e^{i(d-{b'}+2\pi)})$ and $(0,\infty e^{i({a'}-d)})$ does not contain any pole $t=1=e^{2\pi j i}$, where $j\in\ZZ$. So, by applying Cauchy Theorem, one gets from \eqref{equation:EaEbA0} that $A_{-1}=-1$, as being expected in \eqref{equation:EaEb} for $\ell=-1$.

Now, suppose $\ell>-1$, and observe that $\displaystyle\frac{b}{c_\ell}=\frac{c_\ell}{q\,a}\,e^{-2\ell\pi i}$. Similarly, it follows from combining \eqref{equation:EaEbI1a} together with \eqref{equation:EaEbI2} and \eqref{equation:EaEbA} that 
\begin{equation}
 \label{equation:EaEbAl}
A_\ell=-1+\left(\int_0^{\infty e^{i(a'-d)}}-\int_0^{\infty e^{i(d-{b'}-2\ell\pi)}}\right)\frac{1}{1-t}\,\frac{1}{e_q(c_\ell\,t/(q\,a))}\,\frac{dt}{t}\,.
\end{equation}
By applying the Residue Theorem, one obtains from \eqref{equation:EaEbAl} that
$$
A_\ell=-1+2\pi i\sum_{j=1}^{\ell+1}\frac{1}{e_q(c_\ell\,e^{-2j\pi i}/(q\,a))}=-1+2\pi i\sum_{j=1}^{\ell+1}\frac{1}{e_q(c\,e^{(\ell-2j)\pi i}/(q\,a))}\,.
$$
This is exactly what is stated in \eqref{equation:EaEb} combined with  \eqref{equation:EaEbR} for $\ell>0$.

The case $\ell<-1$ can be treated in a similar way, again using \eqref{equation:EaEbAl}. Thus, 
one finishes the proof of Lemma \ref{lemma:EaEb}.  
\end{proof}

\subsection{The $q$-Borel transform of the product of two generalized $q$-Euler series}\label{subsection:qEulerBorel}

In view of Lemma \ref{lemma:Ea}, we shall only consider any product of the form $\hat E_q^{(a,0)}\,\hat E_q^{(b,0)}$, whose differentiations with respect to $a$ and/or to $b$ may lead us to every product  $\hat E_q^{(a,\alpha)}\,\hat E_q^{(b,\beta)}$ for all non-negative integers $\alpha$ and $\beta$.

\begin{lem}\label{lemma:EabEquation}
For any $(a,b)\in\CC^{*2}$, the product series $\hat E_q^{(a,0)}\,\hat E_q^{(b,0)}$ satisfies the following $q$-difference equations:
 \begin{equation}
  \label{equation:EabEquation}
  (x\,\sigma_q-a)(x\,\sigma_q-b)(x^2\,\sigma_q-ab)y=qx^2-ab\,.
 \end{equation}

\end{lem}
\begin{proof}
This comes from considering \eqref{equation:EaEq}. Indeed, letting $\alpha=0$ in \eqref{equation:EaEq} gives that $$(x\,\sigma_q-a)\hat E_q^{(a,0)}(x)=1\,,$$ or, equivalently, $x\,\hat E_q^{(a,0)}(qx)=1+a\,\hat E_q^{(a,0)}(x)$. This implies that 
$$
x^2\,\hat E_q^{(a,0)}(qx)\,\hat E_q^{(b,0)}(qx)=(1+a\,\hat E_q^{(a,0)}(x))(1+b\,\hat E_q^{(b,0)}(x))\,,
$$
which can be put into the following form:
$$(x^2\,\sigma_q-a\,b)\left(\hat E_q^{(a,0)}\,\hat E_q^{(b,0)}\right)(x)=1+a\,\hat E_q^{(a,0)}(x)+b\,\hat E_q^{(b,0)}(x)\,.
$$
Thus, \eqref{equation:EabEquation} is obtained by applying the $q$-difference operator $(x\,\sigma_q-a)(x\sigma_q-b)$ to both sides of the above equation. 
\end{proof}

\begin{lem}\label{lemma:EabBorel}
Given $(a,b)\in\CC^*\times\CC^*$, if $\phi_0(a,b;\xi)=\hat\cB_{q;1}\left(\hat E_q^{(a,0)}\hat E_q^{(b,0)}\right)$, then:
 \begin{equation}\label{equation:phi0}
 \phi_0(a,b;\xi)=-\frac{1}{ab}\,\sum_{n\ge0}\left(1+\frac{a\,q^n}{\xi-a\,q^n}+\frac{b\,q^n}{\xi-b\,q^n}\right)\,q^{-n^2}\,\left(\frac{\xi^{2}}{ab}\right)^{n}\,.
\end{equation}
 Consequently, $\phi_0(a,b;\xi)\in\HH_{q;2}^d$ for any $d\in\RR$ such that $\{a,b\}\cap (0,\infty e^{id})=\emptyset$. 
\end{lem}

\begin{proof} For simplicity, we will write $\phi_0(\xi)$ instead of $\phi_0(a,b;\xi)$. By considering \eqref{equation:qBorelx} for $j=m=1$,
applying the $q$-Borel transform $\hat\cB_{q;1}$ to both sides of \eqref{equation:EabEquation} yields that 
$$
(\xi-a)(\xi-b)(q^{-1}\,\xi^2\,\sigma_q^{-1}-ab)\phi_0(\xi)=\xi^2-ab\,.
$$In other words, 
$\phi_0$
satisfies the following $q$-difference equation:
\begin{equation}\label{equation:CarreBorel}
(1-X)u(\xi) =v(\xi)\,,
\end{equation} 
where $X=\displaystyle\frac{\xi^2}{abq}\sigma_q^{-1}$ and $v(\xi)=\displaystyle\frac{ab-\xi^2}{ab(\xi-a)(\xi-b)}$. By observing that $$\xi^2-ab=(\xi-a)(\xi-b)+a(\xi-b)+b(\xi-a),$$ one obtains that 
\begin{equation}
 \label{equation:Eabv}
 v(\xi)=-\frac{1}{ab}\left(1+\frac{a}{\xi-a}+\frac{b}{\xi-b}\right)\,.
\end{equation}

Now, put  \eqref{equation:CarreBorel} into the form $u=v+Xu$, and iterate this last relation. It follows that
$$
u(\xi)=(1+X+\dots+X^{N-1})v(\xi)+X^Nu(\xi)
$$
for any positive integer $N$. By considering the equality $X^n=\displaystyle\left(\frac{\xi^2}{ab}\right)^nq^{-n^2}\sigma_q^{-n}$ for all $n\in\ZZ$, one deduces from \eqref{equation:Eabv} that 
$$
X^nv(\xi)=-\frac{q^{-n^2}\,\xi^{2n}}{(ab)^{n+1}}\,\left(1+\frac{a\,q^n}{\xi-a\,q^n}+\frac{b\,q^n}{\xi-b\,q^n}\right)\,.
$$In this way, one can put  \eqref{equation:CarreBorel} into the following form: 
$$
u(\xi)=-\sum_{n=0}^{N-1}\frac{q^{-n^2}\,\xi^{2n}}{(ab)^{n+1}}\,\left(1+\frac{a\,q^n}{\xi-a\,q^n}+\frac{b\,q^n}{\xi-b\,q^n}\right)+\frac{q^{-N^2}\xi^{2N}}{(ab)^N}\,u\left(\frac\xi {q^N}\right)\,.
$$
Letting $N\to\infty$, one arrives at the following expression deduced from \eqref{equation:CarreBorel}:
\begin{equation}
 \label{equation:EabN}
u_0(\xi)=-\sum_{n\ge0}\frac{q^{-n^2}\,\xi^{2n}}{(ab)^{n+1}}\,\left(1+\frac{a\,q^n}{\xi-a\,q^n}+\frac{b\,q^n}{\xi-b\,q^n}\right)\,,
\end{equation}
that is well-defined and analytic  in $\CC\setminus\left(a\,q^{\ZZ_{\ge 0}}\cup b\,q^{\ZZ_{\ge 0}}\right)$.

By using the fact that $\phi_0$ is the unique analytic function at $\xi=0$ that satisfies the functional equation in \eqref{equation:CarreBorel}, one finds that $\phi_0=u_0$. This gives \eqref{equation:phi0}, using \eqref{equation:EabN}.

To see the growth at the infinity of $\phi_0$, let $V=V_\epsilon^d$ to be any sector that does not contain neither $a$ nor $b$, and set
\begin{equation}\label{equation:abV}
|a|_V=\min_{\xi\in V}|\xi-a|\,,\quad |b|_V=\min_{\xi\in V}|\xi-b|\,.
\end{equation}
By taking into account the fact that
$$
\left|1+\frac{a\,q^n}{\xi-a\,q^n}+\frac{b\,q^n}{\xi-b\,q^n}\right|\le 1+\frac{|a|}{|a|_V}+\frac{|b|}{|b|_V}
$$for all $\xi\in V$,
one gets from \eqref{equation:phi0} that, as $\xi$ tends to $\infty$ in $V_\epsilon^d$,
\begin{equation}\label{equation:phi0F}
\phi_0(\xi)=O(F(\xi))\quad\hbox{where}\quad F(\xi)=\sum_{n\ge0}q^{-n^2}\,\left(\frac{\xi^{2}}{ab}\right)^n\,.
\end{equation}
Notice that 
$F=\hat\cB_{q;2}\hat f
$, where $\hat f(x)=\displaystyle\sum_{n\ge 0}\left(\frac{x^{2}}{\sqrt q\,ab}\right)^n$.
By applying \cite[Proposition 3.1.4]{zhang1999developpements} for $q^{1/2}$ instead of $q$, one find that $F$ has a $q^{1/2}$-exponential growth of order one at the infinity. Thus, one finishes the proof of Lemma \ref{lemma:EabBorel}.
\end{proof}
 
 \subsection{Towards the summability of the product of two generalized $q$-Euler series}\label{subsection:qEulerLaplace}
 The next step is to consider the second order $q$-Laplace transform of the function $\phi_0(a,b;\xi)$ as follows.
 
 \begin{lem}\label{lemma:phi1d} 
 Let $a$, $b$, $\phi_0(a,b;\xi)$ and $d$ be as in Lemma \ref{lemma:EabBorel}, assume the relation in \eqref{equation:abd}, define  $\zeta_{a,b}^d=q^{1/4}\,\sqrt{a\,b}\,e^{-\pi i}$ with 
 \begin{equation}
  \label{equation:Lambda}
  \Lambda_{a,b}^d=\zeta_{a,b}^d\,e^{\pi i\ZZ}\quad\hbox{and}\quad \Lambda_{a,b}^{d*}=\Lambda_{a,b}^d\setminus\{\zeta_{a,b}^d\}\,,
 \end{equation}
 and set $\phi_1^{d}(a,b;\zeta)=\cL_{q;2}^d(\phi_0(a,b;\xi))(\zeta)$.
 \begin{enumerate}
  \item The function $\phi_1^{d}(a,b;\zeta)$ is well-defined and analytic in $\tilde\CC^*\setminus \Lambda_{a,b}^{d*}$, possessing a simple pole at each point of $\Lambda_{a,b}^{d*}$ and satisfying the following equality for all $\zeta\in\tilde\CC^*\setminus \Lambda_{a,b}^{d*}$:
\begin{equation}
 \label{equation:phi1d}(q^{-1/2}\zeta^{2}-ab)\,\phi_1^{d}(a,b;\zeta) =1+a \,\mathcal{S}_{q;2}^{d}(\hat E_{q^{1/2}}^{(a,0)})(\zeta)+b\,\mathcal{S}_{q;2}^{d}(\hat E_{q^{1/2}}^{(b,0)})(\zeta)\,,
\end{equation}
where $\mathcal{S}_{q;2}^{d}=\cL_{q;2}^d\circ\mathcal{C}^d\circ\hat{\mathcal{B}}_{q;2}=\mathcal{S}_{q^{1/2};1}^{d}$.
\item
Given any $d'\in\RR$ and $\epsilon>0$, one can find $C>0$ such that 
\begin{equation}
 \label{equation:phi1dC}
 |\phi_1^d(a,b;\zeta)|\le \frac C{|\zeta^2-\sqrt q\,a\,b|} 
\end{equation}
for all $\zeta\in V_\epsilon^{d'}$ with $\zeta^2\not=\sqrt q\,a\,b$. 
\item One has $\displaystyle\phi_1^{d}(a,b;\zeta)\in\tilde\HH_{q;2}^{d'}$ for all $\displaystyle d'\in\RR\setminus\left(\frac12\left(\arg(a)+\arg(b)\right)-\pi+\pi\ZZ^*\right)$.
 \end{enumerate}

\end{lem}

\begin{proof} 
(1) Remember that $\phi_0(a,b;\xi)$ is solution to the functional equation stated in \eqref{equation:CarreBorel}, which together with \eqref{equation:Eabv} implies that
\begin{equation}
 \label{equation:Eabphi0}
\left(1-\frac{\xi^2}{a\,b\,q}\sigma_q^{-1}\right)\phi_0(a,b;\xi)=-\frac{1}{a\,b}\left(1+\frac{a}{\xi-a}+\frac{b}{\xi-b}\right)\,.
\end{equation}
From now on, let $q'=q^{1/2}$. By definition, one knows that $\cL_{q;2}^d=\cL_{q';1}^d$. A direct computation shows that 
$$\cL_{q';1}^d(\xi^2\,\sigma_q^{-1})=
\cL_{q';1}^d\left(q'\,(\xi\,\sigma_{q'}^{-1})^2\right)=q'\,\zeta^2\,\cL_{q';1}^d\,.
$$
Thus, applying $\cL_{q;2}^d$ to both side of \eqref{equation:Eabphi0} yields that
\begin{equation}
 \label{equation:phi1dL}
\left(1-\frac{\zeta^2}{a\,b\,q'}\right)\phi_1^d(a,b;\zeta)=-\frac{1}{a\,b}\,\cL_{q';1}^d\left(1+\frac{a}{\xi-a}+\frac{b}{\xi-b}\right)(\zeta)\,.
\end{equation}
Note that $\cL_{q';1}^d(1)=1$.
In view of \eqref{equation:EaS} with $\alpha=0$, it follows from \eqref{equation:phi1dL} that 
$$
\left(1-\frac{\zeta^2}{a\,b\,q'}\right)\phi_1^d(a,b;\zeta)=-\frac{1}{a\,b}\,\left(1+a\,\mathcal{S}_{q';1}^d\hat E_{q'}^{(a,0)}(\zeta)+b\,\mathcal{S}_{q';1}^d\hat E_{q'}^{(b,0)}(\zeta)\right)\,.
$$
Since $\mathcal{S}_{q^{1/2};1}^d=\mathcal{S}_{q;2}^d$, one deduces thus the equality stated in \eqref{equation:phi1d}. 

By noticing that $\zeta^2=\sqrt q\,a\,b$ if and only if $\zeta\in\Lambda_{a,b}^d$ (see \eqref{equation:Lambda}), one deduces from \eqref{equation:phi1d} that $\phi_1^d$ is well-defined and analytic in $\tilde\CC^*\setminus\Lambda_{a,b}^d$,  eventually admitting simple poles in $\Lambda_{a,b}^d$.
Furthermore, if one replaces $q$ with $q^{1/2}$ in Lemma \ref{lemma:EaEb}, one can easily deduce from \eqref{equation:EaEbR} that, given $\ell\in\ZZ$, one has 
$$1+a \,\mathcal{S}_{q;2}^{d}(\hat E_{q^{1/2}}^{(a,0)})(\zeta_{a,b}^d\,e^{\ell\pi i})+b\,\mathcal{S}_{q;2}^{d}(\hat E_{q^{1/2}}^{(b,0)})(\zeta_{a,b}^d\,e^{\ell\pi i})=0
$$
if and only if $\ell=-1$. This implies that $\phi_1^d(a,b;\zeta)$ is well-defined and analytic at $\zeta_{a,b}^d$ and really admits simple poles in the complementary set $\Lambda_{a,b}^{d*}$.

(2)
With regard to the relation expected in \eqref{equation:phi1dC}, one can make use of direct estimates on the $q'$-Laplace transforms appeared in \eqref{equation:phi1dL}, observing that both $Gq$-sums $\mathcal{S}_{q';1}^d(\hat E_{q'}^{(a,0)})$ and $\mathcal{S}_{q';1}^d(\hat E_{q'}^{(b,0)})$ are bounded in any sector $V_{\epsilon}^{d'}$ of the Riemann surface $\tilde\CC^*$ (see \eqref{equation:Eabounded}). 

(3) This comes immediately from combining the relation in \eqref{equation:phi1dC} with the fact that $\phi_1^d(a,b;\zeta)$ is well-defined and analytic at $\zeta_{a,b}^d$.
\end{proof}

\begin{rmk}\label{remark:EaEbphi1}
 By taking into account the condition stated in  \eqref{equation:abd}, one has 
 \begin{equation}
  \label{equation:abdd}
 d\notin \left(\frac12\left(\arg(a)+\arg(b)\right)-\pi+\pi\ZZ^*\right).
 \end{equation}
This together with Lemma \ref{lemma:phi1d} (3) implies that $\phi_1^d(a,b;\zeta)\in\tilde\HH_{q;2}^d$.
\end{rmk}

\subsection{Proof of Theorem \ref{theorem:Eab}}\label{subsection:proofThmEab}

Before starting the proof, we will establish one result about the product of two $q$-exponential functions given as in  \eqref{equation:e_q}.

\begin{lem}\label{lemma:Proofeq} 
Let $e_q$ as in \eqref{equation:e_q}, and let $(\xi_1,\xi_2)\in \tilde\CC^*\times \tilde\CC^*$ and $(\zeta_1,\zeta_2)\in \tilde\CC^*\times \tilde\CC^*$. If
\begin{equation}\label{equation:Proofxi}
\zeta_1=q^{-1/4}\,\sqrt{\frac{\xi_1}{\xi_2}}\,,\quad 
 \zeta_2=q^{1/4}\,\sqrt{\xi_1\xi_2}\,,
\end{equation}
then the following identity holds for any $x\in\tilde\CC^*$:
\begin{equation}\label{equation:Proofeq}
e_q\left(\frac{\xi_1}x\right)\,e_q\left(\frac{\xi_2}x\right)=\frac12\,e_{q^{1/2}}(\zeta_1)\,e_{q^{1/2}}\left(\frac{\zeta_2} x\right)\,.
\end{equation}

\end{lem}

\begin{proof} 
By applying the logarithm application to both sides of \eqref{equation:Proofeq}, considering \eqref{equation:e_q} together with \eqref{equation:Proofxi} implies that \eqref{equation:Proofeq} is equivalent to the following one:
$$
\left(\log\left(\frac{q^{1/2}\,\xi_1}x\right)\right)^2+\left(\log\left(\frac{q^{1/2}\,\xi_2}x\right)\right)^2=\frac{1}{2}\,\left(\log\left(\frac{\xi_1}{\xi_2}\right)\right)^2+\frac12\,\left(\log\left(\frac{q\,\xi_1\,\xi_2}{x^2}\right)\right)^2\,.
$$
This can be obtained by direct computation, noticing that 
$$
(\log A)^2+(\log B)^2=\frac12\,\left(\left(\log\left(\frac AB\right)\right)^2+(\log(A\,B))^2\right)
$$
for all $A$, $B\in\tilde\CC^*$. Thus, one fnishes the proof of Lemma \ref{lemma:Proofeq}.  
\end{proof}

\begin{proof}
[Proof of Theorem \ref{theorem:Eab}] 
Through the whole proof, fix $d\in\RR$ such that $\{a,b\}\cap (0,\infty e^{id}) =\emptyset$,  also assuming the convention stated in \eqref{equation:abd}. 
 
By gathering Lemmas \ref{lemma:EabBorel} and \ref{lemma:phi1d} together with Remark \ref{remark:EaEbphi1}, one finds that the product $\hat E_q^{(a,0)}\,\hat E_q^{(b,0)}$ belongs to the space $\CC\{x\}_{q;(1,2)}^d$. For simplify, we will make use of the following notation:
$
 f^d(c;x)=\mathcal{S}_{q;1}^{d}\hat E_q^{(c,0)}(x)$ for $c=a$ or $b$,  and $
 f^d(a,b;x)=\mathcal{S}_{q;(1,2)}^{d}\hat E_q^{(a,0)}\,\hat E_q^{(b,0)}(x)$.
In this way, the  relation stated in \eqref{equation:EaEbS} of Theorem \ref{theorem:Eab}, that is what we shall prove, can be expressed as follows:
\begin{equation}
 \label{equation:EaEbf}
 f^d(a,b;x)=f^d(a;x)\,f^d(b;x)\,.
\end{equation}
One the one hand, by applying Definition \ref{definition:Gq2} and Lemma \ref{lemma:phi1d}, it follows that, for all $x\in\tilde\CC^*$, 
\begin{equation}\label{equation:Prooffab}
f^d(a,b;x)=\int_0^{\infty e^{id}}\frac{\phi_1^d(a,b;\zeta)}{e_{q^{1/2}}(\zeta/x)}\,\frac{d\zeta}{\zeta}\,,
\end{equation}
where $\phi_1^d(a,b;\zeta)$ is given as in \eqref{equation:phi1d}.
On the other hand, by letting $(a,\alpha)=(a,0)$ or $(b,0)$ in \eqref{equation:EaS}, one can write
\begin{equation}
 \label{equation:Prooffafb} f^d(a;x)\,f^d(b;x)=\int_0^{\infty  e^{\mathbf{i}d}}\!\!\int_0^{\infty e^{\mathbf{i}d}}\frac{1}{(\xi_1-a)(\xi_2-b)}\,\frac1{e_q(\xi_1/x)\,e_q(\xi_2/x)}\,\frac{d\xi_1\,d\xi_2}{\xi_1\,\xi_2}
\end{equation}
for all $x\in\tilde\CC^*$.

Let $\Phi$ : $(\xi_1,\xi_2)\mapsto(\zeta_1,\zeta_2)$ be the homeomorphism from $(0,\infty e^{id})\times(0,\infty e^{id})$ onto $(0,+\infty)\times(0,\infty\,e^{id})$ defined by \eqref{equation:Proofxi}. By noticing 
 $\displaystyle
 \xi_1=\zeta_1\,\zeta_2$ and $\displaystyle \xi_2=\frac{q^{-1/2}\,\zeta_2}{\zeta_1}
 $, the corresponding Jacobian is as follows: 
 $$\left|
 \frac{\partial(\xi_1,\xi_2)}{\partial(\zeta_1,\zeta_2)}\right|=\left|
\begin{array}{ll}
\zeta_2&\zeta_1\cr\displaystyle
  -q^{-1/2}\,\zeta_2/\zeta_1^2& \displaystyle q^{-1/2}/\zeta_1\cr
\end{array}
\right|=2\,q^{-1/2}\,\frac{\zeta_2}{\zeta_1}\,.
 $$
 Thus, one gets that
 $$\frac{d\xi_1\,d\xi_2}{\xi_1\,\xi_2}=2\,q^{-1/2}\,\frac{\zeta_2}{\zeta_1}\,\times\,\frac{d\zeta_1\,d\zeta_2}{\xi_1\,\xi_2}\,=2\,\frac{d\zeta_1\,d\zeta_2}{\zeta_1\,\zeta_2}\,,
 $$what together with \eqref{equation:Proofeq} allows us to transform \eqref{equation:Prooffafb} into the following form:
 \begin{equation}\label{equation:Prooffafb2}
 f^d(a;x)\,f^d(b;x)=\int_0^{\infty\,e^{id}}\!\!\int_0^{+\infty}\frac{\eta(\zeta_1,\zeta_2)}{e_{q^{1/2}}(\zeta_1)\,e_{q^{1/2}}(\zeta_2/x)}\,\frac{d\zeta_1\,d\zeta_2}{\zeta_1\,\zeta_2}\,,
 \end{equation}
 where
 $$\eta(\zeta_1,\zeta_2)=
 \frac{\zeta_1}{(\zeta_1\,\zeta_2-a)(q^{-1/2}\,\zeta_2-b\,\zeta_1)}\,.
 $$
 Therefore, letting 
 \begin{equation}\label{equation:ProofH}
  H(\zeta_2)=\int_0^{+\infty}\frac{\eta(\zeta_1,\zeta_2)}{e_{q^{1/2}}(\zeta_1)}\,\frac{d\zeta_1}{\zeta_1}\,,
 \end{equation}
 the relation in \eqref{equation:Prooffafb2} can be expressed as follows:
  \begin{equation}\label{equation:ProoffafbH}
 f^d(a;x)\,f^d(b;x)=\int_0^{\infty\,e^{id}}\frac{H_2(\zeta_2)}{e_{1/2}(\zeta_2/x)}\,\frac{d\zeta_2}{\zeta_2}\,.
 \end{equation}
  
 Note that the following fractional decomposition with respect to $\zeta_1$ holds:
 \begin{equation}\label{equation:ProofH1}
 \eta(\zeta_1,\zeta_2)=\frac{1}{q^{-1/2}\,\zeta_2^2-a\,b}\,\left(\frac{a}{\zeta_2\,\zeta_1-a}+\frac{\zeta_2}{\zeta_2-q^{1/2}\,b\,\zeta_1}\right)\,.
 \end{equation}
 Letting \eqref{equation:ProofH1} into \eqref{equation:ProofH} yields that
\begin{equation}\label{equation:ProofH2}
  H(\zeta_2)=\frac{1}{q^{-1/2}\,\zeta_2^2-ab}\,\left(I_1(\zeta_2)+I_2(\zeta_2)\right)\,,
 \end{equation}
where
\begin{equation}\label{equation:ProofHI1}
I_1(\zeta_2)=a\,\int_0^{+\infty}\frac{1}{(\zeta_2\,\zeta_1-a)\,e_{q^{1/2}}(\zeta_1)}\,\frac{d\zeta_1}{\zeta_1}
\end{equation}
and 
\begin{equation}\label{equation:ProofHI2}
I_2(\zeta_2)=\int_0^{+\infty}\frac{\zeta_2}{(\zeta_2-q^{1/2}\,b\,\zeta_1)\,e_{q^{1/2}}(\zeta_1)}\,\frac{d\zeta_1}{\zeta_1}\,.
\end{equation}
Putting $\zeta_1=\xi/\zeta_2$ into \eqref{equation:ProofHI1} implies that
\begin{equation}\label{equation:ProofHI1*}
I_1(\zeta_2)=a\,\int_0^{\infty e^{id}}\frac{1}{(\xi-a)\,e_{q^{1/2}}(\xi/\zeta_2)}\,\frac{d\xi}{\xi}=a\,\mathcal{S}_{q^{1/2};1}^d(\hat E_{q^{1/2}}^{(a,0)})(\zeta_2)\,.
\end{equation} 
Besides, by noticing that 
$$ \frac{\zeta_2}{\zeta_2-q^{1/2}\,b\,\zeta_1}=1+\frac{q^{1/2}\,b\,\zeta_1}{\zeta_2-q^{1/2}\,b\,\zeta_1}=1+\frac{b}{q^{-1/2}\,\zeta_2/\zeta_1-b}$$ and that $e_{q^{1/2}}(\zeta_2/(q^{1/2}\,\xi))=e_{q^{1/2}}(\xi/\zeta_2)$, letting $\zeta_1=\zeta_2/(q^{1/2}\,\xi)$ into \eqref{equation:ProofHI2} yields that 
\begin{equation}\label{equation:ProofHI2*}
I_2(\zeta_2)=1+b\,\int_0^{\infty\,e^{id}}\frac{1}{(\xi-b)\,e_{q^{1/2}}(\xi/\zeta_2)}\,\frac{d\xi}{\xi}=1+b\,\mathcal{S}_{q^{1/2};1}^d(\hat E_{q^{1/2}}^{(b,0)})(\zeta_2)\,.
\end{equation}
Hence, it follows from combining \eqref{equation:ProofH2} with  \eqref{equation:ProofHI1*} and \eqref{equation:ProofHI2*} that 
\begin{equation}\label{equation:ProofHS}
H_2(\zeta_2)=\frac{1}{q^{-1/2}\,\zeta_2^2-a\,b}\,\left(1+a\,\mathcal{S}_{q^{1/2};1}^d(\hat E_{q^{1/2}}^{(a,0)})(\zeta_2)+b\,\mathcal{S}_{q^{1/2};1}^d(\hat E_{q^{1/2}}^{(b,0)})(\zeta_2)\right)\,.
\end{equation}
Considering both \eqref{equation:phi1d} and \eqref{equation:ProofHS} implies that $H_2=\phi_1^d$. Finally, comparing \eqref{equation:Prooffab} and \eqref{equation:ProoffafbH} allows us to finishe the proof of Theorem \ref{theorem:Eab}.
\end{proof}

\section{Proof of  Theorem \ref{theorem:intro}}\label{section:proof}

This section aims to proving Theorem \ref{theorem:intro}, by using Theorem \ref{theorem:Eab} together with Theorem \ref{theorem:generic}. For doing that, we will establish a generalization of Theorem \ref{theorem:Eab} for any product power series of the form $\hat E_q^{(a,\alpha)}\,\hat E_q^{(b,\beta)}$; see Theorem \ref{theorem:EabD}.

\subsection{Differentiating with respect to the parameters in the first Borel-plane}\label{subsection:D}

The binomial coefficients $\binom mk$ are usually defined for non-negative integers $m$ and $k$. In what follows, it will be convenient to make use of the following extension:
\begin{equation}\label{equation:Dbinom}
 \binom mk=\left\{ 
\begin{array}{ll}
 1&\hbox{for $m\in\ZZ_{<0}$ and $k=0$,}\cr 
 0&\hbox{for $m\in\ZZ_{<0}$ and $k>0$.} 
\end{array}
 \right.
\end{equation}

\begin{lem}\label{lemma:phi0D}
 Given any $(a,b,\alpha,\beta)\in\CC^*\times\CC^*\times\ZZ_{\ge 0}\times\ZZ_{\ge 0}$, one has
 \begin{equation}\label{equation:phi0D}
  \hat\cB_{q;1}\left(\hat E_q^{(a,\alpha)}\,\hat E_q^{(b,\beta)}\right)(\xi)=\frac{(-1)^{\alpha+\beta}}{a^{\alpha+1}\,b^{\beta+1}}\,\sum_{n\ge 0}C_n(a,b;\xi)\,q^{-n^2}\,\left(\frac{\xi^2}{a\,b}\right)^n\,,
 \end{equation}
 where the sequence $\{C_n(a,b;\xi)\}_{n\ge 0}$ is defined as follows:
 \begin{align}\label{equation:phi0DC}
  C_n(a,b;\xi)=-\binom{n+\alpha}\alpha\,\binom{n+\beta}\beta&+\binom{n+\beta}\beta\,\sum_{\ell=1}^{\alpha+1}\binom{\alpha-\ell+n} {\alpha-\ell+1}\,\left(\frac{a\,q^n}{a\,q^n-\xi}\right)^{\ell}\cr
  &+\binom{n+\alpha}\alpha\,\sum_{\ell=1}^{\beta+1}\binom{\beta-\ell+n} {\beta-\ell+1}\,\left(\frac{b\,q^n}{b\,q^n-\xi}\right)^{\ell}\,.
 \end{align}
 
\end{lem}

\begin{proof} By taking into account Lemma \ref{lemma:EabBorel} and  the fact that the $q$-Borel transform commutes with the derivation for the parameters, considering \eqref{equation:qEulera} yields that
\begin{equation}\label{equation:phi0D00}
 \hat\cB_{q;1}\left(\hat E_q^{(a,\alpha)}\,\hat E_q^{(b,\beta)}\right)(\xi)=\frac{1}{\alpha!\,\beta!}\,\partial_{a,b}^{\alpha,\beta} \phi_0(a,b;\xi)\,.
\end{equation}  
Moreover,
is is easy to see that \eqref{equation:phi0} can be put into the following form:
\begin{equation}\label{equation:phi0(ab)1}
 \phi_0(a,b;\xi)=\sum_{n\ge0}\left(-\frac{1}{(a\,b)^{n+1}}+\frac{q^n}{a^n\,b^{n+1}\,(a\,q^n-\xi)}+\frac{q^n}{a^n\,b^{n+1}(b\,q^n-\xi)}\right)\,q^{-n^2}\,\xi^{2n}\,.
\end{equation}
Since
$$
\partial_{a,b}^{\alpha,\beta}\frac{1}{(a\,b)^{n+1}}=\frac{(-1)^{\alpha+\beta}}{(a\,b)^{n+1}}\,\frac{(n+\alpha)!}{n!\,a^\alpha}\,\frac{(n+\beta)!}{n!\,b^\beta}\,,
$$
one finds that
\begin{equation}\label{equation:phi0D0}
\frac{1}{\alpha!\,\beta!}\,\partial_{a,b}^{\alpha,\beta}\frac{1}{(a\,b)^{n+1}}=\frac{(-1)^{\alpha+\beta}}{a\,b}\,\frac{\binom{n+\alpha}\alpha\,\binom{n+\beta}\beta}{a^\alpha\,b^\beta}\,\frac{1}{(a\,b)^{n}}\,.
\end{equation}
Besides, by noticing that
$$
\partial_{a,b}^{\alpha,\beta}\frac{1}{a^n\,b^{n+1}\,(a\,q^n-\xi)}=\frac{(-1)^{\beta}}{b^{n+1}}\,\frac{(n+\beta)!}{n!\,b^\beta}\,\partial_{a}^{\alpha}\frac{1}{a^n\,(a\,q^n-\xi)}
$$
and
$$
\partial_{a}^{\alpha}\frac{1}{a^n\,(a\,q^n-\xi)}=\frac{(-1)^{\alpha}}{a^n\,(a\,q^n-\xi)}\,\sum_{k=0}^\alpha\binom\alpha k\frac{(n+k-1)!}{(n-1)!\,a^k}\frac{(\alpha-k)!\,q^{(\alpha-k)\,n}}{(a\,q^n-\xi)^{\alpha-k}}\,,
$$
one gets that
\begin{equation}\label{equation:phi0D10}
 \frac{1}{\alpha!\,\beta!}\,\partial_{a,b}^{\alpha,\beta}\frac{1}{a^n\,b^{n+1}\,(a\,q^n-\xi)}=\frac{(-1)^{\alpha+\beta}}{a\,b}\,\frac{\binom{n+\beta}\beta}{b^\beta}\,\sum_{k=0}^\alpha\frac{\binom{n+k-1} {k}\,q^{(\alpha-k)\,n}}{a^{k-1}\,(a\,q^n-\xi)^{\alpha-k+1}}\,\frac{1}{(a\,b)^n}\,.
\end{equation}
Letting $\ell=\alpha-k+1$ into the summation of the right-hand side of \eqref{equation:phi0D10} yields that
\begin{equation}\label{equation:phi0D1}
 \frac{1}{\alpha!\,\beta!}\,\partial_{a,b}^{\alpha,\beta}\frac{1}{a^n\,b^{n+1}\,(a\,q^n-\xi)}=\frac{(-1)^{\alpha+\beta}}{a\,b}\,\frac{\binom{n+\beta}\beta}{b^\beta}\,\sum_{\ell=1}^{\alpha+1}\frac{\binom{\alpha+n-\ell} {\alpha+1-\ell}\,q^{(\ell-1)\,n}}{a^{\alpha-\ell}\,(a\,q^n-\xi)^{\ell}}\,\frac{1}{(a\,b)^n}\,.
\end{equation}
Similarly, one can find that
\begin{equation}\label{equation:phi0D2}
 \frac{1}{\alpha!\,\beta!}\,\partial_{a,b}^{\alpha,\beta}\frac{1}{a^{n+1}\,b^{n}\,(b\,q^n-\xi)}=\frac{(-1)^{\alpha+\beta}}{a\,b}\,\frac{\binom{n+\alpha}\alpha}{a^\alpha}\,\sum_{\ell=1}^{\beta+1}\frac{\binom{\beta+n-\ell} {\beta+1-\ell}\,q^{(\ell-1)\,n}}{b^{\beta-\ell}\,(b\,q^n-\xi)^{\ell}}\,\frac{1}{(a\,b)^n}\,.
\end{equation}
By gathering \eqref{equation:phi0D0}, \eqref{equation:phi0D1} and \eqref{equation:phi0D2}, considering both \eqref{equation:phi0D00} and \eqref{equation:phi0(ab)1} implies \eqref{equation:phi0D}.
\end{proof}

In what follows, we will write
\begin{equation}\label{equation:Dphi0}
 \phi_0^{\alpha,\beta}(a,b;\xi)=\hat\cB_{q;1}\left(\hat E_q^{(a,\alpha)}\,\hat E_q^{(b,\beta)}\right)(\xi)\,.
\end{equation}

\begin{prop}\label{proposition:D1}
 Let $(a,b)\in\CC^*\times\CC^*$, $(\alpha,\beta)\in\ZZ_{\ge 0}\times\ZZ_{\ge 0}$, and let $d\in\RR$ be such that $\{a,b\}\cap(0,\infty e^{id})=\emptyset$. One has $\phi_0^{\alpha,\beta}(a,b;\xi)\in\HH_{q;2}^d$.
\end{prop}

\begin{proof}
 Both $\hat E_q^{(a,\alpha)}$ and $\hat E_q^{(b,\beta)}$ belonging to the $q$-Gevrey space $\CC[[x]]_{q;1}$, the $q$-Borel transform $\phi_0^{\alpha,\beta}(a,b;\xi)$ given in \eqref{equation:Dphi0} represents a germ of analytic function at the origin of the complex plane.  In order to prove that $\phi_0^{\alpha,\beta}(a,b;\xi)\in\HH_{q;2}^d$,  let $\epsilon>0$, suppose $V=V_{\epsilon}^d$ be such that $\{a,b\}\cap V=\emptyset$, and consider the positive constants $|a|_V$ and $|b|_V$ defined  in \eqref{equation:abV}. Set
\begin{equation*}\label{equation:Dab}
 c=\max(|a|,|b|)\,,\quad 
 \gamma=\max(\alpha,\beta)\,,\quad \delta=\min(|a|_V,|b|_V)\,.
\end{equation*} 
By using $\displaystyle \frac{z\,q^n}{z\,q^n-\xi}=\frac{z}{z-q^{-n}\,\xi}$ for $z=a$ or $b$, one can notice that 
\begin{equation}\label{equation:Dabd}
 \max\left(\sup_{\xi\in V}\left|\frac{a\,q^n}{a\,q^n-\xi}\right|,\sup_{\xi\in V}\left|\frac{b\,q^n}{b\,q^n-\xi}\right|\right)\le\frac{c}{\delta}
\end{equation}
for all integer $n$.

In what follows, we will suppose that $\gamma>0$,  the case of $\gamma=0$ being already known with $\alpha=\beta=0$. 
As $\displaystyle\binom{n}k\le n^k$ for any non-negative integers $n$ and $k$,  one deduces from combining \eqref{equation:phi0DC} and \eqref{equation:Dabd} that for any $\xi\in V$, 
  $$
  \left\vert C_n(a,b;\xi)\right\vert\le (n+\gamma)^{2\gamma}+2(n+\gamma)^{2\gamma}\,\sum_{\ell=1}^{\gamma+1}\left(\frac{c}{\delta}\right)^{\ell}\,.
  $$
  This together with the expression given in
\eqref{equation:phi0D} implies that there exists $C>0$ such that
\begin{equation}\label{equation:phi0CV}
\left|\phi_0^{\alpha,\beta}(a,b;\xi)\right|\le C\,\sum_{n\ge 0}(n+\gamma)^{2\gamma}\,q^{-n^2}\,\left(\frac{|\xi|^2}{|a|\,|b|}\right)^n
\end{equation}
for all $\xi\in V$.
So, by comparing \eqref{equation:phi0CV} with \eqref{equation:phi0F}, considering the same argument as what used at the end of the proof of Lemma \ref{lemma:EabBorel} implies that $\phi_0^{\alpha,\beta}(a,b;\xi)\in\HH_{q;2}^d$.
\end{proof}

\subsection{Fractional decomposition arround a pole in the second Borel-plane}\label{subsection:Removable} 
Given $(a,b)\in\CC^*\times\CC^*$ and $(\alpha,\beta)\in\ZZ_{\ge 0}\times\ZZ_{\ge 0}$, we define  $\left\{P_j^{\alpha,\beta}\right\}_{0\le j\le\alpha+\beta}\subset\CC[X,Y]$ as follows: 
\begin{equation}\label{equation:phi1DP}
\partial_{a,b}^{\alpha,\beta} \frac1{q^{-1/2}\zeta^{2}-a\,b}=\frac1{q^{-1/2}\zeta^{2}-a\,b}\,\sum_{j=0}^{\alpha+\beta}\frac{P_j^{\alpha,\beta}(a,b)}{(q^{-1/2}\zeta^{2}-a\,b)^j}\,.
\end{equation}
Il will be convenient to write $P_j^{\alpha,\beta}=0$ when $j>\alpha+\beta$ or $j<0$. A direct computation shows that $P_0^{0,0}=1$, $P_j^{0,0}=0$ for $j>0$, and $P_0^{\alpha,\beta}=0$ for $\alpha+\beta>0$; moreover,
\begin{equation*}\label{equation:phi1DP+1}
P_j^{\alpha+1,\beta}(a,b)=\partial_aP_j^{\alpha,\beta}(a,b)+jbP_{j-1}^{\alpha,\beta}(a,b)\,,\quad
P_j^{\alpha,\beta+1}(a,b)=\partial_bP_j^{\alpha,\beta}(a,b)+jaP_{j-1}^{\alpha,\beta}(a,b)\,.
\end{equation*}

\begin{lem}\label{lemma:phi1D}
Let $\phi_1^{d}(a,b;\zeta)$ be as in Lemma \ref{lemma:phi1d}.
 One has
\begin{equation}\label{equation:phi1DFraction}
 \partial_{a,b}^{\alpha,\beta} \phi_1^{d}(a,b;\zeta) =
\frac1{q^{-1/2}\zeta^{2}-a\,b}\,
\sum_{j=0}^{\alpha+\beta}
\frac{F_j(\zeta)}{(q^{-1/2}\zeta^{2}-a\,b)^j}\,,
\end{equation}
where
\begin{align}\label{equation:phi1DFractionF}
 F_j(\zeta)=P_j^{\alpha,\beta}(a,b)&+\sum_{k=j}^\alpha\frac{\alpha!}{k!}\,\left(k\,P_j^{k-1,\beta}(a,b)+a\,P_j^{k,\beta}(a,b)\right)\,\mathcal{S}_{q;2}^{d}(\hat E_{q^{1/2}}^{(a,\alpha-k)})(\zeta)\nonumber \\
&+\sum_{k=j}^\beta\frac{\beta!}{k!}\,\left(k\,P_j^{\alpha,k-1}(a,b)+b\,P_j^{\alpha,k}(a,b)\right)\,\mathcal{S}_{q;2}^{d}(\hat E_{q^{1/2}}^{(b,\beta-k)})(\zeta)\,.
\end{align}
\end{lem}

\begin{proof} 
By \eqref{equation:phi1d}, it follows that
\begin{equation}\label{equation:phi1(ab)}
\phi_1^{d}(a,b;\zeta) =\frac1{q^{-1/2}\zeta^{2}-ab}\,\left(1+a \,\mathcal{S}_{q;2}^{d}(\hat E_{q^{1/2}}^{(a,0)})(\zeta)+b\,\mathcal{S}_{q;2}^{d}(\hat E_{q^{1/2}}^{(b,0)})(\zeta)\right)\,. 
\end{equation}
We will consider the expression of $\phi_1^d(a,b;\zeta)$ given in \eqref{equation:phi1(ab)}, taking the differential $\partial_{a,b}^{\alpha,\beta}$ of the members of its right-hand side. Since
$$
\partial_{a,b}^{\alpha,\beta}
\frac{a \,\mathcal{S}_{q;2}^{d}(\hat E_{q^{1/2}}^{(a,0)})(\zeta)}{q^{-1/2}\zeta^{2}-a\,b}=\alpha!\,
\partial_{b}^{\beta}\,\sum_{k=0}^\alpha
\frac{1}{k!}\,\partial_a^k\left(\frac{a }{q^{-1/2}\zeta^{2}-a\,b}\right)\,\frac{1}{(\alpha-k)!}\,\partial_a^{\alpha-k}\mathcal{S}_{q;2}^{d}(\hat E_{q^{1/2}}^{(a,0)})(\zeta)\,,
$$
using \eqref{equation:SEadD} gives that
\begin{equation}\label{equation:phi1D1}
\partial_{a,b}^{\alpha,\beta}
\frac{a \,\mathcal{S}_{q;2}^{d}(\hat E_{q^{1/2}}^{(a,0)})(\zeta)}{q^{-1/2}\zeta^{2}-a\,b}=
\sum_{k=0}^\alpha
\frac{\alpha!}{k!}\,\partial_{a,b}^{k,\beta}\left(\frac{a }{q^{-1/2}\zeta^{2}-a\,b}\right)\,\mathcal{S}_{q;2}^{d}(\hat E_{q^{1/2}}^{(a,\alpha-k)})(\zeta)\,.
\end{equation}
Similarly, one has
\begin{equation}\label{equation:phi1D2}
 \partial_{a,b}^{\alpha,\beta} \frac{b \,\mathcal{S}_{q;2}^{d}(\hat E_{q^{1/2}}^{(b,0)})(\zeta)}{q^{-1/2}\zeta^{2}-a\,b}=
\sum_{k=0}^\beta
\frac{\beta!}{k!}\,\partial_{a,b}^{\alpha,k}\left(\frac{b }{q^{-1/2}\zeta^{2}-a\,b}\right)\,\mathcal{S}_{q;2}^{d}(\hat E_{q^{1/2}}^{(b,\beta-k)})(\zeta)\,.
\end{equation}
In view of \eqref{equation:phi1D1} and \eqref{equation:phi1D2}, taking the differential $\partial_{a,b}^{\alpha,\beta}$ on both sides of \eqref{equation:phi1(ab)} yields that
 \begin{align}\label{equation:phi1D}
 \partial_{a,b}^{\alpha,\beta} \phi_1^{d}(a,b;\zeta) =&
\partial_{a,b}^{\alpha,\beta} \frac1{q^{-1/2}\zeta^{2}-a\,b}+
\sum_{k=0}^\alpha
\frac{\alpha!}{k!}\,\partial_{a,b}^{k,\beta}\left(\frac{a }{q^{-1/2}\zeta^{2}-a\,b}\right)\,\mathcal{S}_{q;2}^{d}(\hat E_{q^{1/2}}^{(a,\alpha-k)})(\zeta)\nonumber \\
&+
\sum_{k=0}^\beta
\frac{\beta!}{k!}\,\partial_{a,b}^{\alpha,k}\left(\frac{b }{q^{-1/2}\zeta^{2}-a\,b}\right)\,\mathcal{S}_{q;2}^{d}(\hat E_{q^{1/2}}^{(b,\beta-k)})(\zeta)\,.
 \end{align}
 
 Since
$$
\partial_{a,b}^{k,\beta}\frac{a }{q^{-1/2}\zeta^{2}-a\,b}=k\,\partial_{a,b}^{k-1,\beta}\frac{1 }{q^{-1/2}\zeta^{2}-a\,b}+a\,\partial_{a,b}^{k,\beta}\frac{1 }{q^{-1/2}\zeta^{2}-a\,b}\,,
$$
replacing $(\alpha,\beta)$ with $(k,\beta)$ or $(k-1,\beta)$ in \eqref{equation:phi1DP} yields
\begin{equation}\label{equation:phi1DPk}
\partial_{a,b}^{k,\beta}\frac{a }{q^{-1/2}\zeta^{2}-a\,b}=\frac1{q^{-1/2}\zeta^{2}-a\,b}\,\sum_{j=0}^{k+\beta}\frac{k\,P_j^{k-1,\beta}(a,b)+a\,P_j^{k,\beta}(a,b)}{(q^{-1/2}\zeta^{2}-a\,b)^j}\,.
\end{equation}
Similarly, one has
\begin{equation}\label{equation:phi1DPk1}
\partial_{a,b}^{\alpha,k}\frac{b }{q^{-1/2}\zeta^{2}-a\,b}=\frac1{q^{-1/2}\zeta^{2}-a\,b}\,\sum_{j=0}^{\alpha+k}\frac{b\,P_j^{\alpha,k}(a,b)+k\,P_j^{\alpha,k-1}(a,b)}{(q^{-1/2}\zeta^{2}-a\,b)^j}\,.
\end{equation}
Therefore, putting \eqref{equation:phi1DP}, \eqref{equation:phi1DPk} and \eqref{equation:phi1DPk1} into \eqref{equation:phi1D} yields \eqref{equation:phi1DFraction}.
\end{proof}

As before, let $(a,b,d)\in\CC^*\times\CC^*\times\RR$ be such that $\{a,b\}\cap(0,\infty e^{id})=\emptyset$. By following Proposition \ref{proposition:D1}, if
\begin{equation}\label{equation:Dabpsi}
\phi_1^{d;\alpha,\beta}(a,b;\zeta)=\int_0^{\infty e^{id}}\frac{\phi_0^{\alpha,\beta}(a,b;\xi)}{e_{q^{1/2}}(\xi/\zeta)}\,\frac{d\xi}{\xi}\,,
\end{equation}
then $\phi_1^{d;\alpha,\beta}(a,b;\zeta)$ represents an analytic function for all $\zeta\in\tilde\CC^*$ such that $|\zeta|<\displaystyle q^{1/4}\,\sqrt{\vert a\,b\vert}$.

\begin{thm}\label{theorem:phi1d}
 Let $\phi_1^{d;\alpha,\beta}(a,b;\zeta)$ as in \eqref{equation:Dabpsi}, $d'\in\RR$, and consider the set $\Lambda_{a,b}^d$ given in \eqref{equation:Lambda}. The following assertions hold.
 \begin{enumerate}
  \item If $(0,\infty e^{id'})\cap \Lambda_{a,b}^d=\emptyset$, there exist $\epsilon>0$ and $C>0$ such that 
  \begin{equation}\label{equation:phi1DCV}
\forall \zeta\in V_\epsilon^{d'}\,,\quad \left|\phi_1^{d;\alpha,\beta}(a,b;\zeta)\right|\le\frac{C}{|\zeta^{2}-\sqrt q\,a\,b|}\,.
 \end{equation}
 \item Otherwise, one has $(0,\infty e^{id'})\cap \Lambda_{a,b}^d=\{c\}$ and there exist a finite set  $\{\lambda_0,\dots,\lambda_{\alpha+\beta}\}\subset\CC$ and an analytic function $\psi$ in $V_{\pi/2}^{d'}$ possessing the following properties.
 \begin{enumerate}
\item For any $\zeta\in V_{\pi/2}^{d'}\setminus\{c\}$, one has
  \begin{equation}\label{equation:phi1DCVd}
\phi_1^{d;\alpha,\beta}(a,b;\zeta)=\sum_{j=0}^{\alpha+\beta}\frac{\lambda_j}{(\zeta-c)^{j+1}}+\psi(\zeta)\,.
 \end{equation}
  \item Given any $\epsilon\in(0,\frac\pi2)$,
$\psi$ remains bounded in $V_\epsilon^{d'}$ and one has $\psi(\zeta)=O(\frac{1}{\zeta})$ as $\zeta\to\infty$ in $V_\epsilon^{d'}$.
 \end{enumerate} 
 \end{enumerate}
\end{thm}

\begin{proof}
Since
$\phi_0^{\alpha,\beta}(a,b;\xi)=\partial_{a,b}^{\alpha,\beta}\phi_0(a,b;\xi)$, equation in \eqref{equation:Dabpsi} takes the following form: 
$$
\phi_1^{d;\alpha,\beta}(a,b;\zeta)=\int_0^{\infty e^{id}}\frac{\partial_{a,b}^{\alpha,\beta}\phi_0(a,b;\xi)}{e_{q^{1/2}}(\xi/\zeta)}\,\frac{d\xi}{\xi}\,.
$$ This implies that
\begin{equation}\label{equation:Dabpsi1}
\phi_1^{d;\alpha,\beta}(a,b;\zeta)=\partial_{a,b}^{\alpha,\beta}\phi_1^d(a,b;\zeta)\,,
\end{equation}
what allows one to make use of the expression given in \eqref{equation:phi1DFraction} for $\partial_{a,b}^{\alpha,\beta}\phi_1^d(a,b;\zeta)$.
Furthermore,  
considering the relation stated in \eqref{equation:Eabounded} implies that everyone of the members of both $Gq$-sums families $\left\{\mathcal{S}_{q;2}^{d}(\hat E_{q^{1/2}}^{(a,k)})\right\}_{0\le k\le \alpha}$ and $\left\{\mathcal{S}_{q;2}^{d}(\hat E_{q^{1/2}}^{(b,k)})\right\}_{0\le k\le\beta}$ is bounded over any sector $V_\epsilon^{d'}$ with $\epsilon>0$.
Thus, one obtains from \eqref{equation:phi1DFractionF} that there exists a constant $C_1>0$, depending of $a$, $b$, $\alpha$, $\beta$ and $V_\epsilon^{d'}$, such that
\begin{equation}\label{equation:phi1DPB}
\forall \zeta\in V_\epsilon^{d'}\,,\quad \max_{0\le j\le\alpha+\beta}\left|F_j(\zeta)\right\vert\le C_1\,.
\end{equation}

(1) If $(0,\infty e^{id'})\cap \Lambda_{a,b}^d=\emptyset$, one has $q^{-1/2}\zeta^{2}-a\,b\not=0$ for all $\zeta\in(0,\infty\,e^{id'})$. Let $\epsilon>0$ be enough small in such way that $q^{-1/2}\zeta^{2}-a\,b \not=0$ for all $\zeta\in V_{2\epsilon}^{d'}$, and define $\delta=\inf_{\zeta\in V_\epsilon^{d'}}|\zeta^{2}-\sqrt q\,a\,b|$.
By letting $C=\sqrt q\,(1+\delta+\dots+\delta^{\alpha+\beta})\,C_1$ and taking into account \eqref{equation:phi1DPB}, one deduces 
\eqref{equation:phi1DCV} from \eqref{equation:phi1DFraction}.

(2) By hypothesis, $c$ is the only root of  $q^{-1/2}\zeta^{2}-a\,b=0$ in the half-line $(0,\infty e^{id'})$. One the one hand, gathering \eqref{equation:phi1DPB} together with \eqref{equation:phi1DFraction} implies that, for any $\epsilon>0$, one has
\begin{equation}\label{equation:phiDinfini}
 \phi_1^{d;\alpha,\beta}(a,b;\zeta)=O(\frac{1}{\zeta^2})\,,\quad V_\epsilon^{d'}\ni\zeta\to\infty\,.
\end{equation}
One the other hand, letting $q^{-1/2}\zeta^{2}-a\,b=q^{-1/2}\,(\zeta-c)(\zeta+c)$ into \eqref{equation:phi1DFraction} and considering the Taylor expansion of each function $F_j(\zeta)$ at $\zeta=c$ yields the fractional decomposition as \eqref{equation:phi1DCVd}, where $\psi$ remains analytic and bounded in every relatively-compact subset of $V_{\pi/2}^{d'}$.  Note that $\displaystyle\sum_{j=0}^{\alpha+\beta}\frac{\lambda_j}{(\zeta-c)^{j+1}}= O(\frac{1}{\zeta})$ as $\zeta\to\infty$. With the help of \eqref{equation:phiDinfini}, one finds that $\psi=O(\frac1\zeta)$ as $\zeta\to\infty$. This finishes the proof of Theorem \ref{remark:phi1d}.
\end{proof}

 The expression given in \eqref{equation:phi1DFraction} for the function $\phi_1^{\alpha,\beta}(a,b;\zeta)$ shows that it may have a singularity at any square-root of $\sqrt q\,a\,b$ in $\tilde\CC^*$. The following result says that it remains well-defined and analytic at $\zeta=q^{1/4}\,\sqrt{a\,b}\,e^{-\pi i}$. This generalizes Lemma \ref{lemma:phi1d} (1) in which one had $\alpha=\beta=0$.
 
\begin{rmk}\label{remark:phi1d}
 Let $(a,b,d)\in\CC^*\times \CC^*\times\RR$ to be such that \eqref{equation:abd} holds. Given any couple $(\alpha,\beta)\in\ZZ_{\ge 0}\times\ZZ_{\ge 0}$, the associated function $\phi_1^{d;\alpha,\beta}(a,b;\zeta)$ remains analytic at $\zeta=q^{1/4}\,\sqrt{a\,b}\,e^{-\pi i}$.
\end{rmk}

This will be deduced from Lemma \ref{lemma:noStokes} ;
for more details, see the end of the proof of Theorem \ref{theorem:EabD} given in \S \ref{subsection:EabD} below.

\subsection{One generalization of Theorem \ref{theorem:Eab}}\label{subsection:EabD}
In order to generalize Theorem \ref{theorem:Eab} for $\hat E_q^{(a,\alpha)}\,\hat E_q^{(b,\beta)}$, we shall establish the following result.

\begin{lem}\label{lemma:noStokes}
 Let $\delta\in\RR$, $N\in\ZZ_{>0}$. Let us consider parameters $\{c_1,\dots,c_N\}\subset(0,\infty e^\delta)$, $\{\lambda_1,\dots,\lambda_N\}\subset\CC$,  and $\{\mu_1,\dots,\mu_N\}\subset\ZZ_{\ge0}$. Let $\psi\in\tilde\HH_{q;1}^{\delta}$  and let us set
 \begin{equation}\label{equation:noStokesP}
  \phi(\xi)=\sum_{j=1}^N\frac{\lambda_j}{(\xi-c_j)^{\mu_j+1}}+\psi(\xi)\,.
 \end{equation}
 If $(c_j,\mu_j)\not=(c_k,\mu_k)$ for all $(j,k)$ such that $j\not=k$, then the following conditions are equivalent.
 \begin{enumerate}
  \item $\lambda_j=0$ for $1\le j\le N$.
  \item  $\phi\in\tilde\HH_{q;1}^{\delta}$.
  \item There exist $\delta^-<\delta <\delta^+$ such that $\phi\in\left(\tilde\HH_{q;1}^{\delta^-}\cap \tilde\HH_{q;1}^{\delta^+}\right)$ and, moreover, the identity $\cL_{q;1}^{\delta^-}\phi=\cL_{q;1}^{\delta^+}\phi$ holds.
 \end{enumerate}

\end{lem}

\begin{proof}
 The equivalence $(1)\Leftrightarrow(2)$ and  implication $(2)\Rightarrow(3)$ being trivial, it suffices to prove $(3)\Rightarrow(2)$ or $(3)\Rightarrow(1)$. This can be done by applying Residue Theorem to the contour integral $\left(\cL_{q;1}^{\delta^-}\phi-\cL_{q;1}^{\delta^+}\phi\right)$. Indeed, one deduces from considering \eqref{equation:noStokesP} that
\begin{equation}\label{equation:L+-}
 \cL_{q;1}^{\delta^-}\phi(x)-\cL_{q;1}^{\delta^+}\phi(x)=2\pi i\,\sum_{j=0}^N\frac{\lambda_j}{\mu_j!}\,\left.\partial_\xi^{\mu_j}\left(\frac1{\xi\,e_q(\xi/x)}\right)\right|_{\xi=c_j}\,.
\end{equation}
 Since $\displaystyle\xi\,e_q(\frac\xi x)=\frac{x}{q}\,e_q(\frac{q\xi}{x})$, one can observe that
 $$
 \frac1{\xi\,e_q(\xi/x)}=\frac{q}{x\,\sqrt{2\pi  \log (q)}}\,e^{P(\xi)}\,,\quad\hbox{where}\quad P(\xi)=\frac{(\log(q^{3/2}\,\xi/x))^2}{2\log (q)}\,. 
 $$
If one defines the sequence $\{P_k\}$ by the relation
\begin{equation}\label{equation:eqP}
\partial_\xi^k\left(\frac1{\xi\,e_q(\xi/x)}\right)=\frac{P_k(\xi)}{\xi\,e_q(\xi/x)}\,,
\end{equation}
then
$
P_0=1$ and $P_{k+1}=P'\,P_k+P_k'$ for all $k\ge 0$. One can find that
\begin{equation}\label{equation:eqP1}
 P_k(\xi)=\frac{Q_k(\log(q^{3/2}\,\xi/x))}{(2(\log (q))\,\xi)^k}\,,\quad Q_k\in\CC[X]\,,\quad \deg Q_k=k\,. 
\end{equation}
Thus, the relation in \eqref{equation:L+-} yields that
\begin{equation}\label{equation:L+-1}
 \cL_{q;1}^{\delta^-}\phi(x)-\cL_{q;1}^{\delta^+}\phi(x)=2\pi i\,\sum_{j=0}^N\frac{\lambda_j\,Q_{\mu_j}(\log(q^{3/2}\,c_j/x))}{\mu_j!\,(2(\log (q)))^{\mu_j}\,c_j^{\mu_j+1}\,e_q(c_j/x)}\,\,. 
\end{equation}
Assume $x\to0$, and consider the different asymptotic scales contained in the right-hand side of \eqref{equation:L+-1}. One can easily obtain that le left-hand side is null if and only if $\lambda_j=0$ for all $j$ from $1$ to $N$. This means that $(3)\Rightarrow(1)$.
\end{proof}

\begin{thm}\label{theorem:EabD}
 Let $(a,b,d)\in\CC^*\times\CC^*\times\RR$ be as in Theorem \ref{theorem:Eab}, and let $\alpha$, $\beta\in\ZZ_{\ge 0}$. One has $\hat E_q^{(a,\alpha)}\,\hat E_q^{(b,\beta)}\in\CC\{x\}_{q;(1,2)}^d$ and
 \begin{equation}
  \label{equation:EaEbSD}
  \mathcal{S}_{q;(1,2)}^d(\hat E_q^{(a,\alpha)}\,\hat E_q^{(b,\beta)})=\mathcal{S}_{q;1}^d (\hat E_q^{(a,\alpha)})\,\mathcal{S}_{q;1}^d(\hat E_q^{(b,\beta)})\,.
 \end{equation}
\end{thm}

\begin{proof}
Let $d'\in\RR$ be as in Theorem \ref{theorem:phi1d} (1), and define 
\begin{equation}\label{equation:Dabdf}
 f^{d',d;\alpha,\beta}(a,b;x)=\int_0^{\infty e^{id'}}\frac{\phi_1^{d;\alpha,\beta}(a,b;\zeta)}{e_{q^{1/2}}(\zeta/x)}\,\frac{d\zeta}{\zeta}
\end{equation}
and
\begin{equation}\label{equation:Dabdfd}
 f^{d',d}(a,b;x)=\int_0^{\infty e^{id'}}\frac{\phi_1^{d}(a,b;\zeta)}{e_{q^{1/2}}(\zeta/x)}\,\frac{d\zeta}{\zeta}
\end{equation}
for all $x\in\tilde\CC^*$. 
By taking into account both 
\eqref{equation:phi1DCV} and \eqref{equation:Dabpsi1}, one gets that 
\begin{equation}\label{equation:Dabdf1}
f^{d',d;\alpha,\beta}(a,b;\zeta)=\partial_{a,b}^{\alpha,\beta}f^{d',d}(a,b;x)\,.
\end{equation}

Let $\bar d_{a,b}=\frac12(\arg(a)+\arg(b))$. In view of \eqref{equation:abd}, one has $\bar d_{a,b} \in(d,2\pi)$.
 One remembers that $f^d(a,b;x)$ was defined by means of the ingeral of \eqref{equation:Prooffab}. By following Lemma \ref{lemma:phi1d}, the function $\phi_1^{d}(a,b;\zeta)$ is analytic in the sector $V_{\pi/2}^{\bar d_{a,b}}$. Thus, comparing \eqref{equation:Prooffab} with  \eqref{equation:Dabdf1} implies that 
 \begin{equation}\label{equation:Dabdf2}
   f^{d',d}(a,b;x)=f^{d}(a,b;x)
 \end{equation}
 for any $d'\in\left(\bar d_{a,b}-\frac\pi2,\bar d_{a,b}+\frac\pi2\right)$. Besides,
combining \eqref{equation:Dabdf2} with \eqref{equation:Dabdfd} yields that
\begin{equation}\label{equation:Dabdf3}
  f^{d',d;\alpha,\beta}(a,b;x)=\partial_{a,b}^{\alpha,\beta}f^{d}(a,b;x)
 \end{equation}
 for any $d'\in\left(\bar d_{a,b}-\frac\pi2,\bar d_{a,b}\right)\cup \left(\bar d_{a,b},\bar d_{a,b}+\frac\pi2\right)$. 
 Thus, applying Lemma \ref{lemma:noStokes} to the triplet $\left(q^{1/2},\phi_1^{d;\alpha,\beta}(a,b;\zeta),\bar d_{a,b}\right)$ instead of $(q,\phi,\delta)$ implies that
 $\phi_1^{d;\alpha,\beta}(a,b;\zeta)\in\tilde\HH_{q;2}^{\bar d_{a,b}}$. One gets that $\phi_1^{d;\alpha,\beta}(a,b;\zeta)$ is analytic at $\zeta=q^{1/4}\,\sqrt{a\,b}\,e^{-\pi i}$, as claimed by Remark \ref{remark:phi1d}. As by-product, it follows from combining this with Theorem \ref{theorem:phi1d} that $\hat E_q^{(a,\alpha)}\,\hat E_q^{(b,\beta)}\in\CC\{x\}_{q;(1,2)}^d$. At the same time, letting $d'=d$ into \eqref{equation:Dabdf1} and using \eqref{equation:EaEbf} yileds \eqref{equation:EaEbSD}, what permits us to achieve the proof of Theorem \ref{theorem:EabD}.
\end{proof}

\subsection{End of the proof of Theorem \ref{theorem:intro}}\label{subsection:proofMainThm} Let $\hat f_{1}$ and $\hat f_{2}$ be two generic $q$-Gevrey series. In view of the notational convention given in \eqref{equation:Ea<}, one knows that $\hat E_q^{(a,-1)}=1$ for any $a\in\CC^*$. Thus, by applyinig \eqref{equation:decomposition} to $\hat f_1$ and $\hat f_2$, one can find  $m$, $n\in\ZZ_{\ge0}$, $g_0, \dots, g_m$, $h_0,\dots, h_n\in\CC\{x\}$, $a_1, \dots, a_m$, $b_1, \dots, b_n\in\CC^*$ and $ \alpha_1, \dots, \alpha_m$, $\beta_1, \dots, \beta_n\in\ZZ_{\ge 0}$ such that
 $$
   \hat f_1=\sum_{j=0}^mg_j\,\hat E_q^{(a_j,\alpha_j)}\,,\quad 
      \hat f_2=\sum_{k=0}^nh_k\,\hat E_q^{(b_k,\beta_k)}\,,
   $$
   where $\alpha_0=\beta_0=-1$ and where $a_0$ and $b_0$ may be arbitrarily chosen in $\CC^*$.
 This implies that 
 \begin{equation}\label{equation:thmprooff12}
\hat f_1\,\hat f_2= \sum_{j=0}^m\sum_{k=0}^ng_j\,h_k\,\hat E_q^{(a_j,\alpha_j)}\,\hat E_q^{(b_k,\beta_k)}\,.
 \end{equation}
 Set $
 \mathbb{S}=\{\arg(a_j)\in[0,2\pi):1\le j\le m\}\cup \{\arg(b_k)\in[0,2\pi):1\le k\le n\}$.
 By taking into account the fact that $\CC\{x\}_{q;(1,2)}^d$ constitutes a $\CC\{x\}$-module (see Propositions \ref{prop:Gqsummable} \& \ref{prop:2Gqsummable}), applying Theorem \ref{theorem:EabD} to each term of the right-hand side of \eqref{equation:thmprooff12} yields that $\hat f_1\,\hat f_2$ is $Gq$-summable of order $(1,2)$ in any direction $d\notin\mathbb{S}$ $\bmod$ $2\pi$.
 Finally, gathering \eqref{equation:thmprooff12} together with \eqref{equation:EaEbSD} allows one to obtain the eqality expected in
 \eqref{equation:f1f2}, with $\mathbb{S}=\{\delta_1, \dots, \delta_M\}$. In this way, we finish the proof of Theorem \ref{theorem:intro}.

 \bibliographystyle{amsalpha}
\bibliography{biblio}
\end{document}